\documentclass{amsart}
\usepackage[english]{babel}
\usepackage{graphicx}
\usepackage{amssymb}
\usepackage{amsmath}
\usepackage[foot]{amsaddr}
\usepackage{amsfonts}
\usepackage{bbm}
\usepackage{xcolor}
\usepackage{float}
\usepackage{algorithm}
\usepackage{algpseudocode}

\usepackage{amsmath}
\usepackage{amssymb}
\usepackage{mathrsfs}
\usepackage{color}
\usepackage{graphicx}
\usepackage[percent]{overpic}
\usepackage{url}
\usepackage{enumerate}
\usepackage{overpic}
\usepackage{amsthm}
\usepackage{moreverb}
\usepackage[pdftex,colorlinks,bookmarksopen,bookmarksnumbered,citecolor=red,urlcolor=red]{hyperref}

\def\0{\boldsymbol{0}}

\def\cl {\nonumber \\}
\def\el {\nonumber }

\newtheorem{rem}{Remark}[section]

\newcommand{\bm}[1]{\mbox{\boldmath{$#1$}}}
\def\div{\nabla\cdot}

\usepackage{tikz}
\usepackage[top=1in, bottom=1.25in, left=1.25in, right=1.25in]{geometry}

\graphicspath{{images/}}

\begin{document}

\title[A linear filter regularization for POD-based ROMs of the QGE]{A linear filter regularization for POD-based Reduced Order Models of the quasi-geostrophic equations
}

\author{Michele Girfoglio$^1$, Annalisa Quaini$^2$ and Gianluigi Rozza$^1$}
\address{$^1$ mathLab, Mathematics Area, SISSA, via Bonomea 265, I-34136 Trieste, Italy}
\address{$^2$ Department of Mathematics, University of Houston, Houston TX 77204, USA}

\begin{abstract}
We propose a regularization for Reduced Order Models (ROMs) of the quasi-geostro\-phic equations (QGE) to increase accuracy when the Proper Orthogonal Decomposition (POD) modes retained to construct the reduced basis are insufficient to describe the system dynamics. Our regularization is based on the so-called BV-$\alpha$ model, which modifies the nonlinear term in the QGE and adds a linear differential filter for the vorticity.
To show the effectiveness of the BV-$\alpha$ model for ROM closure, we compare the results computed by a POD-Galerkin ROM with and without regularization for the classical double-gyre wind forcing benchmark. 
Our numerical results show that the solution computed by the regularized ROM is more accurate, even when the retained POD modes account for a small percentage of the eigenvalue energy. Additionally, we show that, although computationally more expensive that the ROM with no regularization, the regularized ROM 
is still a competitive alternative to full order simulations of the QGE.
\end{abstract}

\maketitle


\vskip .5cm
\centerline{Dedicated to the memory of Roland Glowinski}
\vskip .5cm

\section{Introduction}\label{sec:intro}
During his long and distinguished career, Roland Glowinski has given outstanding contributions to the development of several methodologies (e.g., nonlinear least squares methods, domain decomposition methods, and fictitious domain methods) with applications to a wide range of problems. One of the fields he has extensively contributed to throughout his career is Computational Fluid Dynamics. Among his most cited works in this field, there is a handbook on finite element methods for incompressible viscous flow \cite{glowinski2003finite}, which is a great example of clear, precise, and unambiguous scientific prose. For some of his works on the Stokes and Navier-Stokes equations, Glowinski has chosen the stream function-vorticity formulation \cite{Glowinski76,Glowinski87,Glowinski88,DEAN1991117,Glowinski91}. This paper focuses on the stream function-potential vorticity formulation of the quasi-geostrophic equations, which are
a simplification of the Navier-Stokes equations used for ocean modeling. 

Ocean flows are characterized by the evolution 
of flow structures (eddies and vortices) with a broad range of spatial scales, the larger scales being of the order of hundreds or thousands kilometers 
and the smaller scales less than 1 mm big. This poses a serious challenge at the computational level, especially in certain flow regimes. Two non-dimensional numbers are often used to describe the ocean flow regime: the Reynolds number $Re$ and the Rossby number $Ro$. The Reynolds number is the ratio of inertial forces to viscous forces, while 
the Rossby number weighs the inertial force
over the Coriolis force. Ocean flows with large $Re$
and small $Ro$ are particularly challenging as they require very fine computational meshes to resolve all the eddy scales, leading to a prohibitive computational cost.

In order to contain the computational cost of ocean flow simulations, assumptions are introduced at the level of the the model describing the dynamics. One simplified model is given by the shallow water equations, which are obtained from the Navier–Stokes equations under the assumption that the horizontal length scale for the problem is much
greater than its vertical length scale. Typically satisfied by ocean flows on large domains since maximum ocean depth is about 10 km, this assumption allows to  average the Navier–Stokes equations over the depth and get rid of the vertical dimension. In the limit of small $Ro$, i.e., when the inertial forces are negligible with respect to the Coriolis and pressure forces, one can further simplify the shallow water equations to obtain the quasi-geostrophic equations (QGE). The name for this model comes from the fact that for $Ro = 0$ one recovers geostrophic flow. See, e.g., \cite{Vallis2006, Cushman-Roisin2011, McWilliams2006} for mathematical and physical fundamentals, \cite{San2012, Carere2021, Strazzullo2017} for some advanced applications and \cite{QGE_review} for a recent review on the QGE.

Although the QGE represent one of the simplest models for geophysical flows, their numerical simulation is still non-trivial. In fact, when the 
the Munk scale (a length that depends on $Ro/Re$) is small, computational simulations require very fine meshes because the mesh size has to be smaller than the Munk scale. Since often long time intervals have to be simulated, the overall computational cost becomes prohibitive. Hence, the need for techniques that reduce the computational cost of simulations for small Rossby numbers 
remains pressing. 

Reduced order models (ROMs) are inexpensive surrogates for expensive models that are built based on a 
relatively few solutions of the latter model and for which the expense
incurred in the construction process is then amortized over many solutions of the surrogate \cite{ModelOrderReduction2, ModelOrderReduction, Rozza2008}. ROMs have been applied to more efficiently treat problems in uncertainty quantification, control and optimization, and other settings that require multiple simulations, 
with applications ranging from biomedical to naval engineering. Among all the existing approaches for ROM development, the proper orthogonal decomposition (POD) is one of the most successful. POD allows to extract the dominant modes from a database of 
high-fidelity numerical solutions. Such modes are used to construct a reduced basis. Then, a way to build the reduced model is by
projecting the governing equations onto the space spanned by the reduced basis. See, e.g., \cite{ModelOrderReduction2, ModelOrderReduction, Rozza2008}. In this paper, we propose a POD-based ROM for the GQE when time is the only parameter.

It is well known that the number of POD basis functions is usually small (meaning, $\mathcal{O}$(10)) for computational problems dominated by diffusion, i.e., for flows characterized by small $Re$. Since the size of the reduced problem depends on the size of the POD basis, POD-based ROMs are very efficient surrogate models for low $Re$ flows. However, realistic geophysical flows are
dominated by convection (i.e., $Re$ is high) and thus
the dimension of the POD basis is expected to
be large. Since the use of large POD basis implies limited computational savings or none at all, one is forced to work with an insufficient number of POD modes in order to indeed reduce the computational time. However, when the number of modes is not enough to capture the relevant flow dynamics, the ROM solution is unphysical, typically displaying spurious numerical oscillations. One way to fix this issue is the introduction of a closure model. 

Closure models, also called regularizations, aim at capturing the effect of the truncated POD modes. In the literature, one can find several strategies to obtain ROM closure models for the QGE. Among these strategies, we mention large eddy simulation (LES) \cite{Selten1995, San2014, Rahman2019, Mou2020, Mou2021, Xie2018bis, Wang2012}, machine learning \cite{Rahman2018, Rahman2019bis, San2018ML}, and stochastic mode reduction \cite{Franzke2005}.
In this paper, we propose a novel ROM closure of the LES type and compare it with the corresponding ROM with no closure to stress the importance of regularization. The particular LES approach that we follow is called BV-$\alpha$, from the fact that the QGE are also known as barotropic vorticity (BV) equations and parameter $\alpha$ (the filtering radius) is of critical importance as we shall see.

The BV-$\alpha$ model has been widely used as a replacement of the QGE, i.e., as
a full order model \cite{Nadiga2001, Holm2003, Monteiro2015, Monteiro2014, Girfoglio1}. By modifying the nonlinear term in the QGE and adding a differential filter for the vorticity (which can be linear \cite{Nadiga2001, Holm2003, Monteiro2015, Monteiro2014} or nonlinear \cite{Girfoglio1}), the BV-$\alpha$ model circumvents the need for a mesh size smaller than the Munk scale. Thanks to this, the BV-$\alpha$ model provides a physical computed solution with coarser meshes than needed by the QGE. Here, we adopt the linear BV-$\alpha$ model for ROM closure. To the best of our knowledge, 
this idea has not been investigated so far.

There are a few other differences between our regularized ROM approach and previous works. The high fidelity solutions for the construction of the reduced basis are obtained by direct numerical simulations of the QGE with an efficient Finite Volume method \cite{Girfoglio1, Girfoglio2}. Other authors have chosen Finite Element methods
\cite{Monteiro2015, Monteiro2014} or Finite Difference methods \cite{Nadiga2001, Holm2003}.  We consider the formulation of the QGE in terms of potential vorticity, instead of standard vorticity as in all previous work.
This choice is justified by the fact that the potential vorticity satisfies a conservation equation which can be exactly enforced by our Finite Volume method
at the discrete level. In addition, we consider different coefficients for the ROM approximation of the potential vorticity and stream function. This leads to two important consequences. First, the stream function basis functions do not depend on
the particular vorticity basis functions. Instead, they are computed directly from the stream function
high-fidelity solutions. Second, the reduced spaces for the stream function and
vorticity can have different dimensions.

For the assesement of the proposed ROM approach, we
consider the classical double-gyre wind forcing benchmark \textcolor{red}{[13, 14, 15, 16, 21,
22]}. We present numerical results for two cases with the same Munk scale: i) Rossby number $Ro = 0.0036$,
Reynolds number $Re = 450$ and ii) Rossby number $Ro = 0.008$, Reynolds number $Re = 1000$. The higher $Re$ in the second case makes the computation of the high-fidelity solutions more challenging, while the smaller $Ro$ of the first case introduces difficulties at the ROM level
as we will show in Sec.~\ref{sec:results}.

The rest of the paper is organized as follows. In Sec.~\ref{sec:FOM}, we introduce the QGE and discuss their time and space discretization. Sec.~\ref{sec:ROM} presents the POD-based ROM and the new closure based on the linear BV-$\alpha$ model.
Numerical results are reported in Sec.~\ref{sec:results}. Conclusions are drawn in Sec.~\ref{sec:End}, where we also presents some future
perspectives.
\section{Mathematical Models}\label{sec:models}


In order to state the quasi-geostrophic equations, let $\Omega$ be a fixed two-dimensional spatial domain and ($t_0$, $T$) a time interval of interest. Let $\omega$ be the standard vorticity (i.e., the curl of the velocity field) and $q = Ro~\omega + y$ the potential vorticity, where $y$ is the non-dimensional vertical coordinate. In addition, we denote with $\psi$ the stream function and set $\bm{\psi} = (0, 0, \psi)$. Then, the QGE in non-dimensional variables read: find $\psi$ and $q$ such that
\begin{align}
\frac{\partial q}{\partial t}+ \div \left(\left(\nabla \times \bm{\psi}\right) q \right) - \dfrac{1}{Re} \Delta q  = F \quad \mbox{ in }\Omega \times (t_0,T), \label{eq:BVE1}\\
-Ro \Delta \psi + y = q   \quad \mbox{ in }\Omega \times (t_0,T), \label{eq:BVE2}
\end{align}
where $Re$ is the Reynolds number, and  $F$ denotes an external forcing. In the rest of the paper, we will focus on external forces that depend exclusively on space. 

Problem \eqref{eq:BVE1}-\eqref{eq:BVE2} needs to be supplemented with proper boundary and initial conditions. We impose
\begin{align}
\psi &= 0 \quad \quad \quad \mbox{ on }\partial \Omega \times (t_0,T), \label{eq:BV5_comp} \\
q &= y \quad \quad \quad \mbox{ on }\partial \Omega \times (t_0,T), \label{eq:BV5_comp2} \\
q(x,y,t_0) &= q_0 = y ~~ \mbox{ in } \Omega, \label{eq:BV5_comp3}
\end{align}
which are a rather standard choice \cite{Nadiga2001,Holm2003,Monteiro2015,Monteiro2014,San2012, Girfoglio1}. Notice that \eqref{eq:BV5_comp}-\eqref{eq:BV5_comp3} are equivalent to $\psi = \omega = 0$ 
on $\partial \Omega$ and $\omega(x,y,t_0) = 0$ in $\Omega$.

The Direct Numerical Simulation (DNS) of the QGE requires a mesh with mesh size $h$ smaller than the Munk scale:
\begin{align}
\delta_M = L \, \sqrt[3]{\dfrac{Ro}{Re}}, \label{eq:munk}
\end{align}
where $L$ is a characteristic length. When a mesh with $h < \delta_M$ cannot be afforded because the associated computational cost would be prohibitive or simply too high, one needs to find a way to model the effects of the scales lost to mesh under-refinement. A possible way to do that is to couple the QGE with a differential filter. The resulting model, called BV-\emph{$\alpha$} \cite{Nadiga2001, Holm2003, Monteiro2015, Monteiro2014}, reads:
find $\psi$, $q$, and $\overline{q}$ such that
\begin{align}
\frac{\partial q}{\partial t} + \div \left(\left(\nabla \times \bm{\psi}\right) q \right) - \dfrac{1}{Re} \Delta q &= F \quad \mbox{ in }\Omega \times (t_0,T), \label{eq:BV1_comp11}\\
-\alpha^2\Delta \overline{q} +\overline{q} &= q  \quad ~~ {\rm in}~\Omega \times
(t_0,T), \label{eq:BV2_comp33} \\
-Ro \Delta \psi + y &= \overline{q} \quad ~ \mbox{ in }\Omega \times (t_0,T), \label{eq:BV2_comp22}
\end{align}
where $\overline{q}$ is the \emph{filtered vorticity} and $\alpha$ can be interpreted as a \emph{filtering radius} (i.e., the radius of the neighborhood where the filter extracts information from the resolved scales). 
The differential filter leverages an elliptic
operator that acts as a spatial filter by damping the spurious and nonphysical oscillations exhibited by the numerical solution on coarse meshes. The price that one has to pay to obtain a physical solution on a coarse grid is the addition of one equation, i.e., eq.~\eqref{eq:BV2_comp33}.

We supplement problem \eqref{eq:BV1_comp11}-\eqref{eq:BV2_comp22} with initial data \eqref{eq:BV5_comp3} and boundary conditions boundary conditions \eqref{eq:BV5_comp}-\eqref{eq:BV5_comp2} plus the additional boundary condition 
\begin{align}
\overline{q} &= y \quad \quad \quad \mbox{ on }\partial \Omega \times (t_0,T). \label{eq:BC_qbar}
\end{align}

While model \eqref{eq:BV1_comp11}-\eqref{eq:BV2_comp22} represents certainly an improvement over model \eqref{eq:BVE1}-\eqref{eq:BVE2} when coarse meshes are used, its effectivity remains limited for severely under-refined meshes. A better alternative is the nonlinear BV-$\alpha$ model introduced in \cite{Girfoglio1}. However, we do not consider the nonlinear BV-$\alpha$ model in this manuscript since nonlinearities pose extra challanges at the ROM level.

\section{The Full Order Method}\label{sec:FOM}
Although we will use both the QGE \eqref{eq:BVE1}-\eqref{eq:BVE2} and the 
BV-$\alpha$ model \eqref{eq:BV1_comp11}-\eqref{eq:BV2_comp22}
to devise the ROMs, at the FOM level we consider only the QGE, i.e., our full order model is given by eqs.~\eqref{eq:BVE1}-\eqref{eq:BVE2} endowed with boundary conditions \eqref{eq:BV5_comp}-\eqref{eq:BV5_comp2} and initial data \eqref{eq:BV5_comp3}. This section is devoted to the time and space discretization of our full order model.

Let us start with the time discretization of eq.~\eqref{eq:BVE1}-\eqref{eq:BVE2}. Let $\Delta t \in \mathbb{R}$, $t^n = t_0 + n \Delta t$, with $n = 0, ..., N_T$ and $T = t_0 + N_T \Delta t$. We denote by $f^n$ the approximation of a generic quantity $f$ at the time $t^n$. 
Problem~\eqref{eq:BVE1}-\eqref{eq:BVE2} discretized in time by a Backward Differentiation Formula of order 1 (BDF1) reads: given $q^0= q_0$, for $n \geq 0$ find the solution $(\psi^{n+1}, q^{n+1})$ to
\begin{align}
\frac{1}{\Delta t}\, q^{n+1} + \div \left(\left(\nabla \times \bm{\psi}^{n+1}\right)q^{n+1}\right) - \dfrac{1}{Re}\Delta q^{n+1} & = b^{n+1},\label{eq:disc_filter_ns-1} \\
-Ro \nabla \psi^{n+1} + y &= q^{n+1},\label{eq:disc_filter_ns-2}
\end{align}
where $b^{n+1} = F + q^n/\Delta t$. 

In order to contain the computational cost required
to approximate the solution to problem \eqref{eq:disc_filter_ns-1}-\eqref{eq:disc_filter_ns-2}, we opt for a segregated algorithm. 
Given the potential vorticity $q^{n}$, at $t^{n+1}$ such algorithm requires to:

\begin{itemize}
    \item [i)] Find the potential vorticity $q^{n+1}$ such that
    \begin{align}
\frac{1}{\Delta t}\, q^{n+1} + \div \left(\left(\nabla \times \bm{\psi}^{n}\right)q^{n+1}\right) - \dfrac{1}{Re}\Delta q^{n+1} & = b^{n+1},\label{eq:disc_filter_ns-1-bis}
     \end{align}
where $\bm{\psi}^{n+1}$ in \eqref{eq:disc_filter_ns-1} is replaced by a linear extrapolation, i.e. 
$\bm{\psi}^{n}$.
    \item [ii)] Find the stream function $\psi^{n+1}$ such that
    \begin{align}
-Ro \nabla \psi^{n+1} + y = q^{n+1}.\label{eq:disc_filter_ns-2-bis}
\end{align}
\end{itemize}

\begin{rem}
The results in Sec.~\ref{sec:results} have been obtained with the BDF1 scheme. For this reason, the algorithm is presented with this choice of temporal discretization. Other schemes are possible. For example, in \cite{Girfoglio1} we reported results given by the BDF2 scheme as well. Therein, we noticed that while BDF1 smooths certain magnitude peaks, the results for BDF1 and BDF2 agree in terms of pattern formation, average kinetic energy, and amplitude of kinetic energy oscillations. 
\end{rem}

For the space discretization of problem  (\ref{eq:disc_filter_ns-1-bis})-(\ref{eq:disc_filter_ns-2-bis}), we adopt a Finite Volume (FV) approximation that is derived directly from the integral form of the governing equations. For this purpose, 
we partition the computational domain $\Omega$ into cells or control volumes $\Omega_i$,
with $i = 1, \dots, N_{c}$, where $N_{c}$ is the total number of cells in the mesh. Let  \textbf{A}$_j$ be the surface vector of each face of the control volume, 
with $j = 1, \dots, M$.
Then, the discretized form of eq.~\eqref{eq:disc_filter_ns-1-bis}, divided by the control volume 
$\Omega_i$, can be written as:
\begin{align}\label{eq:evolveFV-1.1_disc}
\frac{1}{\Delta t}\, q^{n+1}_i &+ \sum_j^{} \varphi^n_j q^{n+1}_{i,j} - \dfrac{1}{Re} \sum_j^{} (\nabla q^{n+1}_i)_j \cdot \textbf{A}_j  = b^{n+1}_i, \quad \varphi_j^{n} = \left(\nabla \times \bm{\psi}_j^{n}\right) \cdot \textbf{A}_j,
\end{align}
where $q_i^{n+1}$ and $b_i^{n+1}$ are the average potential vorticity and source term in control volume $\Omega_i$ and $q_{i,j}^{n+1}$ the potential vorticity associated to the centroid of face $j$ normalized by the
volume of $\Omega_i$.
The discretized form of eq.~\eqref{eq:disc_filter_ns-2-bis}  divided by the control volume 
$\Omega_i$ reads:
\begin{align}
  - Ro \sum_j \left(\nabla\psi_i^{n+1}\right)_j \cdot \textbf{A}_j + y_i = q_i^{n+1} \label{eq:QGE3},
\end{align}
with $\psi_i^{n+1}$ denoting the average stream function in control volume $\Omega_i$ and $y_i$ is the vertical coordinate of the centriod. For further details, we refer the reader to \cite{Girfoglio1, Girfoglio2}.

For the implementation of the numerical scheme described in this section, we chose the finite volume C++ library OpenFOAM\textsuperscript{\textregistered} \cite{Weller1998}.


\section{The Reduced Order Models}\label{sec:ROM}

We assume that any solution to problem \eqref{eq:BVE1}-\eqref{eq:BVE2} can be approximated as a linear combination of a ``small'' number of global basis
functions dependent on space variables only, multiplied by scalar coefficients that depend on time and/or other parameters, which can be physical or geometrical. As mentioned in Sec.~\ref{sec:intro}, in this paper we are interested in the time reconstruction of the flow field
with no other parameter. Hence, the solution $(q, \psi)$ to eq.~\eqref{eq:BVE1}-\eqref{eq:BVE2} is approximated by the reduced solution $(q_r, \psi_r)$, with
\begin{align}
q_r = \sum_{i=1}^{N_{q}^r} \beta_i(t) {\varphi}_i(\bm{x}), \quad 
\psi_r = \sum_{i=1}^{N_{\psi}^r} \gamma_i(t) \xi_i(\bm{x}). \label{eq:ROM_1} 
\end{align}
In \eqref{eq:ROM_1}, $N_{\Phi}^r$ denotes the cardinality of a reduced basis for the space field $\Phi = \{q, \psi\}$. We remark that we consider different coefficients for the approximation of the potential vorticity $q$ and stream function $\psi$ fields. This is unlike the previous works 
where the same coefficients are used for both variables, i.e., 
$\beta_i(t) = \gamma_i(t)$. If one uses the BV-$\alpha$ model \eqref{eq:BV1_comp11}-\eqref{eq:BV2_comp22}, the solution $(q, \psi, \overline{q})$ is approximated by the reduced solution $(q_r, \psi_r, \overline{q}_r)$ with $(q_r, \psi_r)$ as
in \eqref{eq:ROM_1} and
\begin{align}
\overline{q}_r = \sum_{i=1}^{N_{\overline{q}}^r} \overline{\beta_i}(t) {\varphi}_i(\bm{x}). \label{eq:ROM_2} 
\end{align}

Extending our methodology to include phyisical parameters (e.g., the Reynolds number) is rather straightforward as such parameters can be handled in the same way we handle time. On the other hand, geometrical parameters require a different treatment \cite{ModelOrderReduction2, ModelOrderReduction, Rozza2008}. 

\subsection{The POD algorithm}
There exist several techniques
in the literature to generate the reduced basis spaces. Some examples are Proper Orthogonal Decomposition (POD), the Proper Generalized Decomposition, and the Reduced Basis with a greedy sampling strategy.
See, e.g., \cite{ModelOrderReduction,ChinestaEnc2017,Chinesta2011, Dumon20111387,Kalashnikova_ROMcomprohtua,quarteroniRB2016,Rozza2008}. 
We generate the reduced basis spaces with the method of snapshots, briefly described hereafter. 


Let $(q_h, \psi_h)$ be the solution computed with the FOM described in Sec.~\ref{sec:FOM}. Eq.~\eqref{eq:evolveFV-1.1_disc}-\eqref{eq:QGE3} get solved 
at every time step, however not all solutions are retained as snapshots. Indeed, only the solutions at time instant $t^j \in \{t^1, \dots, t^{N_t}\} \subset (t_0, T]$, with $N_T$ a multiple of $N_t$, are stored into the snapshot matrices:
\begin{align}\label{eq:space}
\bm{\mathcal{S}}_{{{\Phi}}} = [{{\Phi}}( t^1), \dots, {{\Phi}}(t^{N_t})] \in \mathbb{R}^{N_{\Phi}^h \times N_t} \quad
\text{for} \quad {{\Phi}} = \{q_h, \psi_h\},
\end{align}
where 
$N_{\Phi}^h$
is the dimension of the full-order space $\Phi$ belong to. 
The POD problem consists in finding, for each value of the dimension of the POD space $N_{POD} = 1, \dots, N_t$, the scalar coefficients $a_1^1, \dots, a_1^{N_t}, \dots, a_{N_t}^1, \dots, a_{N_t}^{N_t}$ and functions ${\zeta}_1, \dots, {\zeta}_{N_t}$ that minimize the error between the snapshots and their projection onto the POD basis. In the $L^2$-norm, such problem reads
\begin{align}
 \text{arg min} \sum_{i=1}^{N_t} \left\lVert {{\Phi}_i} - \sum_{k=1}^{N_{POD}} a_i^k {\zeta}_k \right\rVert \quad \forall N_{POD} = 1, \dots, N_t    \cl
\text{with} \quad ({\zeta}_i, {\zeta}_j)_{L_2(\Omega)} = \delta_{i,j} \quad \forall i,j = 1, \dots, N_t. \label{eq:min_prob}
\end{align}

It can be shown \cite{Kunisch2002492} that problem~\eqref{eq:min_prob} is equivalent to the following eigenvalue problem
\begin{align}
\bm{\mathcal{C}}^{{\Phi}} \bm{Q}^{{\Phi}} &= \bm{Q}^{{\Phi}} \bm{\Lambda}^{{\Phi}}, \label{eq:eigen_prob} \\
\mathcal{C}_{ij}^\Phi &= ({\Phi}_i, {\Phi}_j)_{L_2(\Omega)} \quad \text{for} \quad i,j = 1, \dots, N_t,
\end{align}
where $\bm{\mathcal{C}}^{{\Phi}}$ is the correlation matrix computed from the snapshot matrix $\bm{\mathcal{S}}_{{{\Phi}}}$, $\bm{Q}^{{\Phi}}$ is the matrix of eigenvectors and $\bm{\Lambda}^{{\Phi}}$ is a diagonal matrix whose diagonal entries are the
eigenvalues of $\bm{\mathcal{C}}^{{\Phi}}$. 
Then, the basis functions are obtained as follows:
\begin{align}\label{eq:basis_func}
{\zeta}_i = \dfrac{1}{N_t \Lambda_i^\Phi} \sum_{j=1}^{N_t} {\Phi}_j Q_{ij}^\Phi.
\end{align}
The resulting POD modes are:
\begin{align}\label{eq:spaces}
L_\Phi = [{\zeta}_1, \dots, {\zeta}_{N_\Phi^r}] \in \mathbb{R}^{N_\Phi^h \times N_\Phi^r}.
\end{align}
The values of $N_\Phi^r < N_t$ are chosen according to reach a user-provided threshold $\varepsilon_\Phi$ for the  cumulative energy of the eigenvalues associated to field $\Phi$:
\begin{equation}
\frac{\sum_{i=1}^{N^r_\Phi}\Lambda_i^\Phi}{\sum_{i=1}^{N_t}\Lambda_i^\Phi} \ge \varepsilon_\Phi.
\label{eq:energy}
\end{equation}

In the following, we will consider two approaches for the reduced order model, which share the same offline stage (i.e., collection of snapshots from the QGE with the method described in Sec.~\ref{sec:FOM} and POD procedure described above) but they differ at the online stage. We call these approaches
\begin{enumerate}
\item QGE-QGE ROM: the system to be solved online results from Galerkin projection of the QGE on the reduced (POD) space;
\item QGE-BV-$\alpha$ ROM: Galerkin projection of the BV-$\alpha$ model on the POD space provides the system that has to be solved online.
\end{enumerate}
The QGE-BV-$\alpha$ ROM is the regularized ROM we introduce in this paper and the QGE-QGE ROM is its non-regularized counterpart. 

\subsection{QGE-QGE ROM}\label{sec:QGE-QGE}
By projecting the QGE onto the reduced space, we obtain the
following system: find ($\psi_r^{n+1}$, $q_r^{n+1}$) that solves
\begin{align}
\left(\frac{1}{\Delta t}\, q_r^{n+1} + \div \left(\left(\nabla \times \bm{\psi}_r^{n}\right)q_r^{n+1}\right) - \dfrac{1}{Re}\Delta q_r^{n+1} - b_r^{n+1}, \varphi_i \right)_{L_2(\Omega)} & = 0, \quad i = 1, \dots, N^r_q, \label{eq:red-1.1} \\
\left(-Ro \Delta \psi_r^{n+1} + y - q_r^{n+1}, \xi_i \right)_{L_2(\Omega)} &= 0, \quad i = 1, \dots, N^r_\psi, \label{eq:red-1.2}
\end{align}
where $b_r^{n+1} = F_r + q_r^n/\Delta t$ 
and ${\varphi}_i$ and $\xi_i$ are the basis functions in \eqref{eq:ROM_1}. 
During the online phase of the QGE-QGE ROM, at time $t^{n+1}$
system \eqref{eq:red-1.1}-\eqref{eq:red-1.2}  has to be solved.

In order to write the algebraic system associated with the reduced problem \eqref{eq:red-1.1}-\eqref{eq:red-1.2}, we
introduce the following matrices:
\begin{align}
&M_{r_{ij}} = ({\varphi}_i, {\varphi}_j)_{L_2(\Omega)}, \quad \widetilde{M}_{r_{ij}} = ({\xi}_i, {\varphi}_j)_{L_2(\Omega)}, \quad A_{r_{ij}} = ({\varphi}_i, \Delta {\varphi}_j)_{L_2(\Omega)}, \label{eq:matrices_evolve1} \\
&B_{r_{ij}} = ({\xi}_i, \Delta \xi_j)_{L_2(\Omega)}, \quad G_{r_{ijk}} = (\varphi_i, \div \left(\left(\nabla \times \xi_j\right) \varphi_k\right))_{L_2(\Omega)}, \quad Y_{r_{ij}} = ({\varphi}_i, y_j)_{L_2(\Omega)}. \label{eq:matrices_evolve2} 
\end{align} 
Then, the matrix form of eq.~\eqref{eq:red-1.1}-\eqref{eq:red-1.2} reads: given $\bm{\beta}^{n}$ and $\bm{\gamma}^{n}$ 
find vectors $\bm{\beta}^{n+1}$ and $\bm{\gamma}^{n+1}$, i.e., the vectors whose entries are the values of coefficients $\beta_i$ and $\gamma_i$ in \eqref{eq:ROM_1} at time $t^{n+1}$, such that
\begin{align}
&\bm{M}_r\left(\dfrac{\bm{\beta}^{n+1} - \bm{\beta}^{n}}{\Delta t}\right) + \left(\bm{\gamma}^n\right)^T\bm{G}_r\bm{\beta}^{n+1} - \dfrac{1}{Re} \bm{A}_r \bm{\beta}^{n+1} = \bm{h}, \label{eq:reduced_1} \\
& - Ro \bm{B}_r \bm{\gamma}^{n+1} + \bm{Y}_r - \widetilde{\bm{M}_r} \bm{\beta}^{n+1} = 0,
\label{eq:reduced_2}
\end{align}
where the entries of vector $\bm{h}$ are $h_i = ({\varphi}_i, F)_{L_2(\Omega)}$.


For clarity, Algorithm \ref{alg:1} presents the pseudocode for QGE-QGE ROM. Lines 2-6 describe the offline stage, while lines 7-9 are perfomed online.

\begin{algorithm}
\caption{Pseudocode for QGE-QGE ROM}\label{alg:01}
\begin{algorithmic}[1]
\State{$q_0, N_{T}$}\Comment{Inputs needed}
\For{$n \in \{0, \dots, N_T - 1\}$}\Comment{Time loop}
\State{Solve system \eqref{eq:evolveFV-1.1_disc}-\eqref{eq:QGE3} } \Comment{QGE simulation}
\EndFor
\State{$\{q_i\}_{i=1}^{N_t} \subseteq \{q^k\}_{k=1}^{N_{T}}$ 
\quad $\{\psi_i\}_{i=1}^{N_t} \subseteq \{\psi^k\}_{k=1}^{N_{T}}$} \Comment{Snapshot collection}
\State{$\mathbb Q^r \doteq$ \text{POD}($\{q_i\}_{i=1}^{N_t}$)} \quad $\Psi^r \doteq$ \text{POD}($\{\psi_i\}_{i=1}^{N_t}$) 
\Comment{POD for vorticity and stream function spaces}
\For{$n \in \{0, \dots, N_T - 1\}$}\Comment{Time loop}
\State{Solve system \eqref{eq:reduced_1}-\eqref{eq:reduced_2}} \Comment{
Standard Galerkin projection}
\EndFor
\end{algorithmic}
\label{alg:1}
\end{algorithm}

\subsection{QGE-BV-$\alpha$ ROM}\label{sec:QGE_BV}
Galerkin projection of the BV-$\alpha$ model onto the reduced space gives us the
following system: find ($\psi_r^{n+1}$, $q_r^{n+1},\overline{q}_r^{n+1})$ that solves

\begin{align}
\left(\frac{1}{\Delta t}\, q_r^{n+1} + \div \left(\left(\nabla \times \bm{\psi}_r^{n}\right)q_r^{n+1}\right) - \dfrac{1}{Re}\Delta q_r^{n+1} - b_r^{n+1}, \varphi_i \right)_{L_2(\Omega)} & = 0, \quad i = 1, \dots, N^r_q,\label{eq:red-1.1_1} \\
\left(-\alpha^2 \Delta \overline{q}_r^{n+1} + \overline{q}_r^{n+1} - q_r^{n+1}, \varphi_i \right)_{L_2(\Omega)} & = 0, \quad i = 1, \dots, N^r_{q}, \label{eq:red-1.2_2} \\
\left(-Ro \Delta \psi_r^{n+1} + y - q_r^{n+1}, \xi_i \right)_{L_2(\Omega)} &= 0, \quad i = 1, \dots, N^r_\psi. \label{eq:red-1.2_3}
\end{align}

Using the matrices defined in \eqref{eq:matrices_evolve1}-\eqref{eq:matrices_evolve2}, we can write the matrix of problem \eqref{eq:red-1.1_1}-\eqref{eq:red-1.2_3} as follows: given $\bm{\beta}^{n}$ and $\bm{\gamma}^{n}$ 
find vectors $\bm{\beta}^{n+1}$ and $\bm{\gamma}^{n+1}$ such that
\begin{align}
\bm{M}_r\left(\dfrac{\bm{\beta}^{n+1} - \bm{\beta}^{n}}{\Delta t}\right) + \left(\bm{\gamma}^n\right)^T\bm{G}_r\bm{\beta}^{n+1} - \dfrac{1}{Re} \bm{A}_r \bm{\beta}^{n+1} & = \bm{h}, \label{eq:reduced_11} \\
-\alpha^2 \bm{A}_r \overline{\bm{\beta}}^{n+1} +  \bm{M}_r \left(\overline{\bm{\beta}}^{n+1} - \bm{\beta}^{n+1} \right) & = 0 \label{eq:reduced_33} \\
-Ro \bm{B}_r \bm{\gamma}^{n+1} + \bm{Y}_r  -\widetilde{\bm{M}_r} \overline{\bm{\beta}}^{n+1} & = 0, \label{eq:reduced_22}
\end{align}
where  vector $\bm{h}$ is the same introduced in Sec.~\ref{sec:QGE-QGE}.

Algorithm \ref{alg:2} presents the pseudocode for QGE-BV-$\alpha$ ROM. Notice that the lines 2-6 (offline phase) are the same as in Algorithm \ref{alg:1}, while
lines 7-9 (online phase) are different.

\begin{algorithm}
\caption{Pseudocode for QGE-BV-$\alpha$ ROM}\label{alg:01}
\begin{algorithmic}[1]
\State{$q_0, N_{T}$}\Comment{Inputs needed}
\For{$n \in \{0, \dots, N_T - 1\}$}\Comment{Time loop}
\State{Solve system \eqref{eq:evolveFV-1.1_disc}-\eqref{eq:QGE3} } \Comment{QGE simulation}
\EndFor
\State{$\{q_i\}_{i=1}^{N_t} \subseteq \{q^k\}_{k=1}^{N_{T}}$ 
\quad $\{\psi_i\}_{i=1}^{N_t} \subseteq \{\psi^k\}_{k=1}^{N_{T}}$} \Comment{Snapshot collection}
\State{$\mathbb Q^r \doteq$ \text{POD}($\{q_i\}_{i=1}^{N_t}$)\quad $\Psi^r \doteq$ \text{POD}($\{\psi_i\}_{i=1}^{N_t}$)} 
\Comment{POD for vorticity and stream function spaces}
\For{$n \in \{0, \dots, N_T - 1\}$}\Comment{Time loop}
\State{Solve system \eqref{eq:reduced_11}-\eqref{eq:reduced_22}} \Comment{
BV-$\alpha$ at the reduced level}
\EndFor
\end{algorithmic}
\label{alg:2}
\end{algorithm}

The QGE-BV-$\alpha$ ROM uses the differential filter
at reduced order level as an eddy viscosity closure approach to stabilize the resulting surrogate model. 
The underlying  analogy is the relationship 
between LES and truncated modal projection \cite{Wang2012}.
While filter regularization for ROMs has been studied for toy models like Burger's equation \cite{Iliescu2016} and more complex models like the incompressible Navier-Stokes equations \cite{Wells2017, Xie2016, Gunzburger2019b, Gunzburger2019, Girfoglio_ROM_fluids, Girfoglio_JCP, Strazzullo2021}, to the best of our knowledge it is the first time that it is proposed for the QGE equations.

The major difference when the BV-$\alpha$ model is used for regularization vs when it is used as FOM
is setting of $\alpha$. For the BV-$\alpha$ model as FOM one takes $\alpha \sim h$ 
\cite{Nadiga2001, Holm2003, Monteiro2015, Monteiro2014}, 
while the optimal value of $\alpha$ for regularization is determined by requiring that the FOM solution and regularized ROM solution are as close as possible in average \cite{Wells2017}. See Remark \ref{rem_alpha}.

\subsection{Treatment of the Dirichlet boundary conditions}
In order to homogeneize the snapshots for $q$ and thus make them independent of the
boundary conditions, we use
a classical approach called lifting function method \cite{Star2019, Girfoglio_JCP}. The idea is to modify the vorticity snapshots as follows: instead of taking the computed $q_h$ snapshots, one takes
\begin{align}
q'_h = q_h - \widetilde{q}_h
\end{align}
where $\widetilde{q}_h$ is the time average of the $q_h$ snapshots, called lifting function. The POD is applied to the $q'_h$ snapshots, i.e., the vorticity snapshots satisfying the homogeneous boundary conditions. Then, the lifting function is added back to the reduced vorticity $q_r$ and reduced filtered vorticity $\overline{q}_r$:
\begin{align}
q_r = \widetilde{q}_h + \sum_{i=1}^{N_{q}^r} \beta_i(t) {\varphi}_i(\bm{x}), \quad
\overline{q}_r = \widetilde{q}_h + \sum_{i=1}^{N_{q}^r} \overline{\beta_i}(t) {\varphi}_i(\bm{x}). \el
\end{align}

\section{Numerical results}\label{sec:results}

In order to assess and compare the QGE-QGE and QGE-BV-$\alpha$ ROMs, we consider the well-known double-gyre wind forcing benchmark, which has been widely studied both at full and reduced order level \cite{Nadiga2001, Holm2003, Greathbatch2000, San2011, Monteiro2015, Monteiro2014}.
The computational domain is the 
$[0, 1] \times [-1, 1]$ rectangle and the forcing term is set to $F = \sin(\pi y)$. 
The time interval of interest is [$t_0$, $T$] = [10, 80].
We will focus on two test cases that have the same Munk scale \eqref{eq:munk}:
\begin{itemize}
\item[-] Case 1: $Ro = 0.0036$ and $Re = 450$;
\item[-] Case 2: $Ro = 0.008$ and $Re = 1000$.
\end{itemize}
While Case 1 has been used to test several ROM closure models \cite{Rahman2019bis, Rahman2018, Rahman2019, Mou2020, QGE_review, Mou2021, Xie2018bis, Wang2012},
to the best of our knowledge it is the first time that a larger value of $Re$ as in Case 2
is considered for the same purpose. 
The interest in Case 2 lies in understanding the performance of our ROM approaches when the same Munk scale stays the same (i.e., same $Ro$ to $Re$ ratio) but the Kolmogorov scale \cite{Kolmogorov41-1, Kolmogorov41-2} decreases
(i.e., smaller $Re$). 
For a study of Case 2 at the full order level, see \cite{Girfoglio1}.


The quantities of interest that we will use to 
test the performance of the ROM approaches are
the kinetic energy of the system $E$:
\begin{equation}
E = \dfrac{1}{2}\int_\Omega  \left( \left(\dfrac{\partial \psi}{\partial y}\right)^2 + \left(\dfrac{\partial \psi}{\partial x}\right)^2 \right)d\Omega
\end{equation}\label{eq:kin_energy}
and the relative $L^2$-norm error $\varepsilon$: 
\begin{equation}\label{eq:normerror}
\varepsilon = \dfrac{||\widetilde{\psi}_{FOM} - \widetilde{\psi}_{ROM}||_{L^2(\Omega)}}{||\widetilde{\psi}_{FOM}||_{L^2(\Omega)}}
\end{equation}
of the time-averaged stream function $\widetilde{\psi}$, which is defined as:
\begin{equation}\label{eq:time_psi}
\widetilde{\psi} = \dfrac{1}{T - t_0}\int_{t_0}^{T} \psi dt.
\end{equation} 
We focus on the stream function, and not on vorticity, because it is $\widetilde{\psi}$
that displays the well-known four gyre structure
that is hard to capture when the numerical method is not accurate. For a study of the potential vorticity in this benchmark, we refer to \cite{Girfoglio1}.

\begin{rem}\label{rem_alpha}
As mentioned in Sec.~\ref{sec:QGE_BV}, the optimal value of $\alpha$ varies whether the BV-$\alpha$ model is used as FOM or ROM. 
We set $\alpha$ for the QGE-BV-$\alpha$ ROM by trial and error, i.e., we try several values and choose the one that minimizes error $\varepsilon$ defined in \eqref{eq:normerror}.
Although computationally cheap, this procedure, which also used in \cite{Wells2017}
for incompressible flow problems, 
could be improved by, e.g., a heuristic criterion.
\end{rem}


Following \cite{Mou2020, San2014, QGE_review, San2011}, we collect 700 FOM snapshots, i.e., one every 0.1 time unit, for the training of the ROM in both cases. The FOM snapshots are computed using structured, orthogonal meshes with a level of refinement that will be specified for each case.
 

\subsection{Results for Case 1}\label{sec:case1}


For a rigorous DNS, one should use a mesh size $h<\delta_M = 0.02$. Indeed, in \cite{San2011} the finite element solutions are computed with 
a 256 $\times$ 512 mesh, which satisfies the condition on the mesh size. In \cite{Girfoglio1}, 
we successfully reproduced the results from \cite{San2011} using the same mesh and a FV approach. In addition, we showed that our FV approach provides accurate approximated solutions of the QGE model with a 16 $\times$ 32 mesh, although it does not satisfy the condition on $h$. We speculate that the reason why we obtain accurate solutions on meshes coarser than necessary is the exact conservation yielded by FV methods at the discrete level. 
Thus, in order to reduce the offline cost without sacrificing accuracy we collect the FOM snapshots with the 16 $\times$ 32 mesh. We set $\Delta t = 1e-4$. 

\begin{figure}[htb]
\centering
\begin{overpic}[width=0.5\textwidth]
    {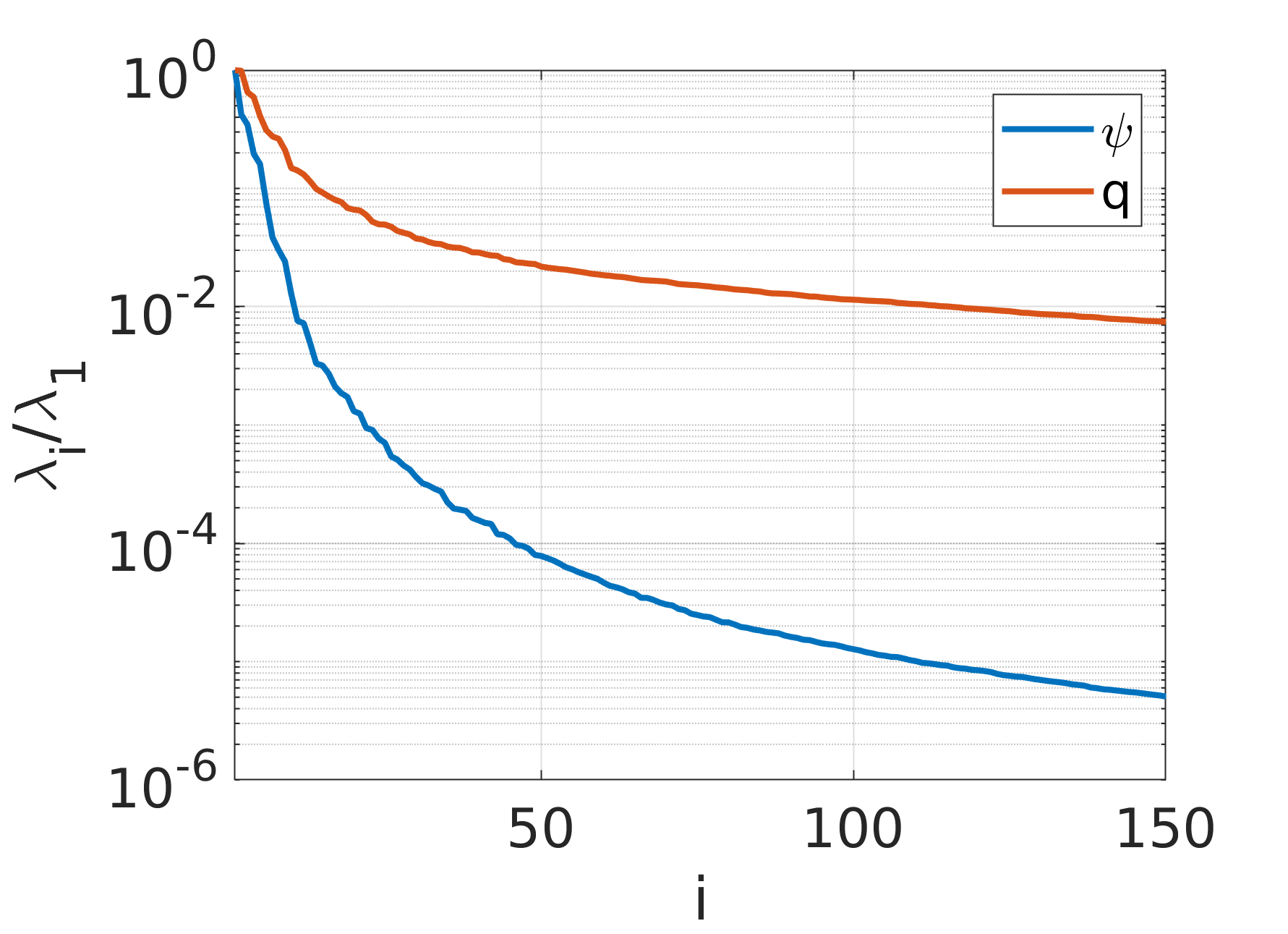}
    \end{overpic}
    \caption{Case 1: eigenvalue decay for stream function and vorticity.}
    \label{fig:eig1}
\end{figure}

We start by displaying the eigenvalue decay for stream function and the  potential vorticity in 
Figure \ref{fig:eig1}. 
Notice how much slower such decay is for $q$ than for $\psi$. This difference in the rapidity of the eigenvalue decay for the two variables makes this problem challenging for ROMs. 
Indeed, with $N^r_\psi = 10$ (i.e., 10 modes for $\psi$) one captures 98\% of the eigenvalue energy \eqref{eq:energy}, while with $N^r_q = 30$ (i.e., 30 modes for $q$) one retains only 70\% of energy \eqref{eq:energy}. 
This is in line with 
\cite{QGE_review}, which reports that 
30 modes for $\omega$ (standard vorticity, instead of potential vorticity) are needed to capture
78\% of the eigenvalue energy. Comparison with the number of modes taken in other works, e.g. \cite{Rahman2019bis, Rahman2018, Rahman2019, Mou2020, Mou2021}, is hindered not only because we state the problem in $q$ rather than $\omega$, but also because we use stream function basis functions that are independent from the vorticity basis functions. Typically, one retains 99.99\% 
of the eigenvalue energy for each variable. For this problem, that would mean $N^r_\psi = 200$ and $N^r_q = 502$, which are too large to lead to any meaningful reduction of the computational time. Thus, we will work with smaller values of  $N^r_\psi$ and $N^r_q$, as done in all previous works.

\begin{figure}[htb]
    \centering
    \begin{overpic}[width=0.2\textwidth]{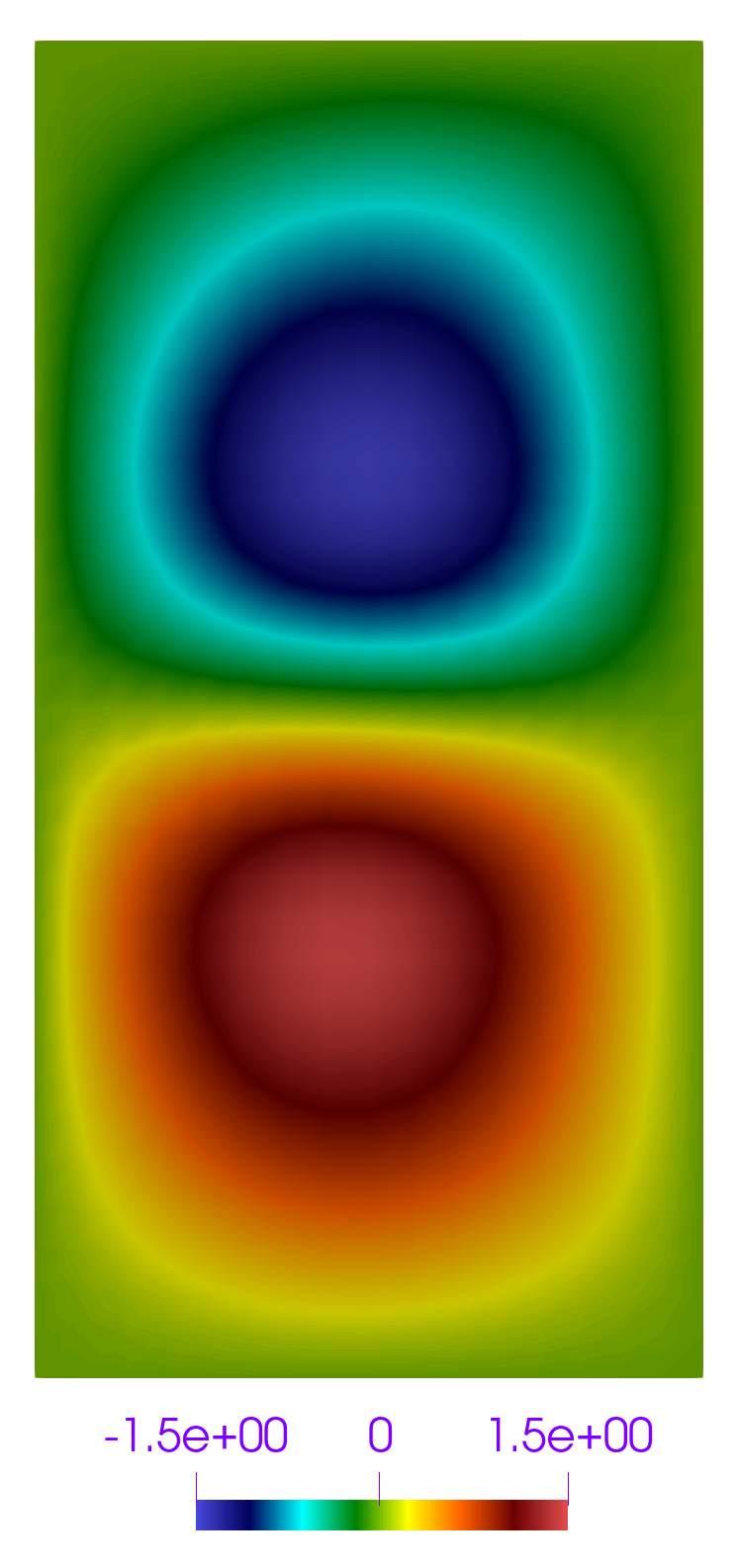}
    \end{overpic}
\put(-48,190){$\xi_1$}
    \begin{overpic}[width=0.2\textwidth]{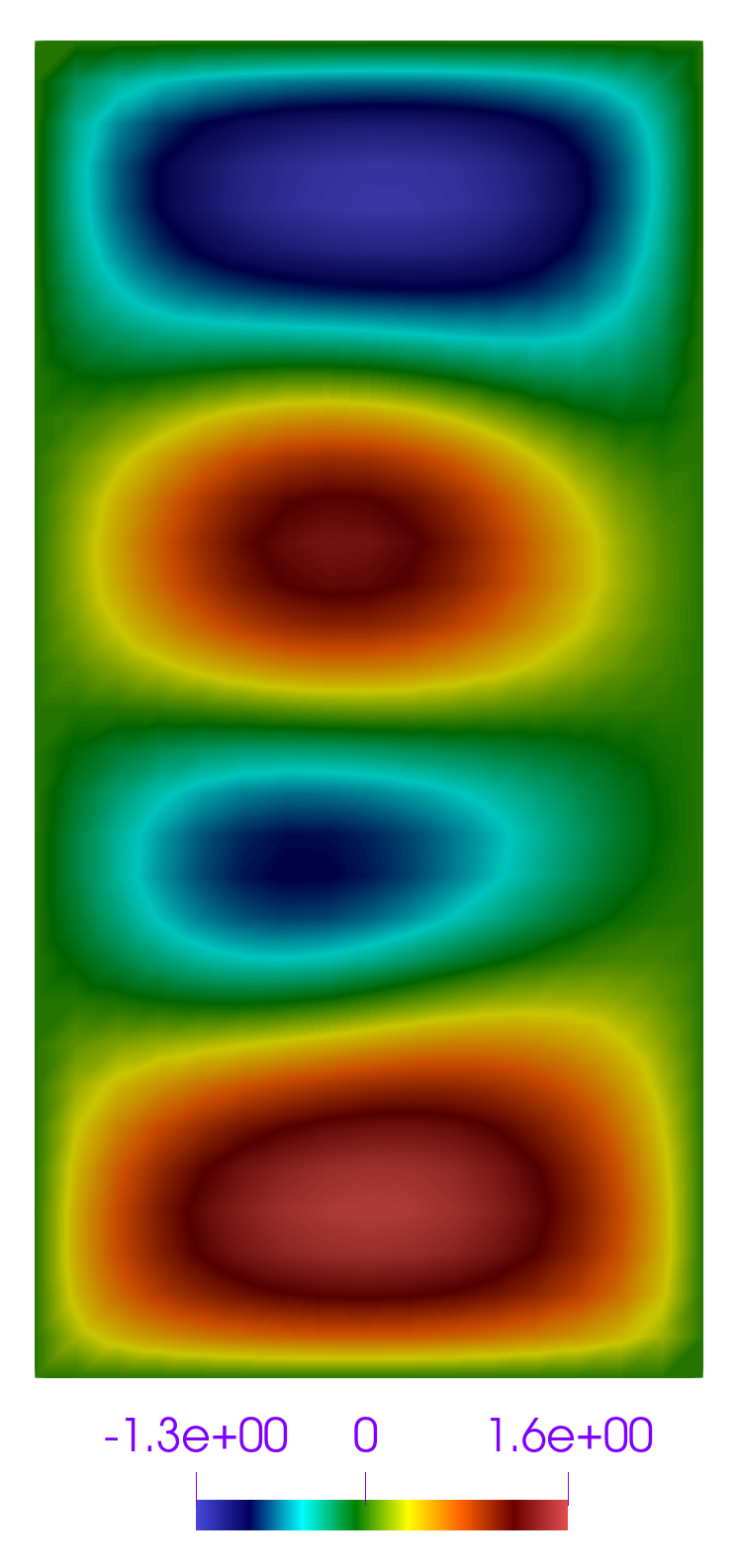}
    \end{overpic}
\put(-48,190){$\xi_3$}
    \begin{overpic}[width=0.2\textwidth]{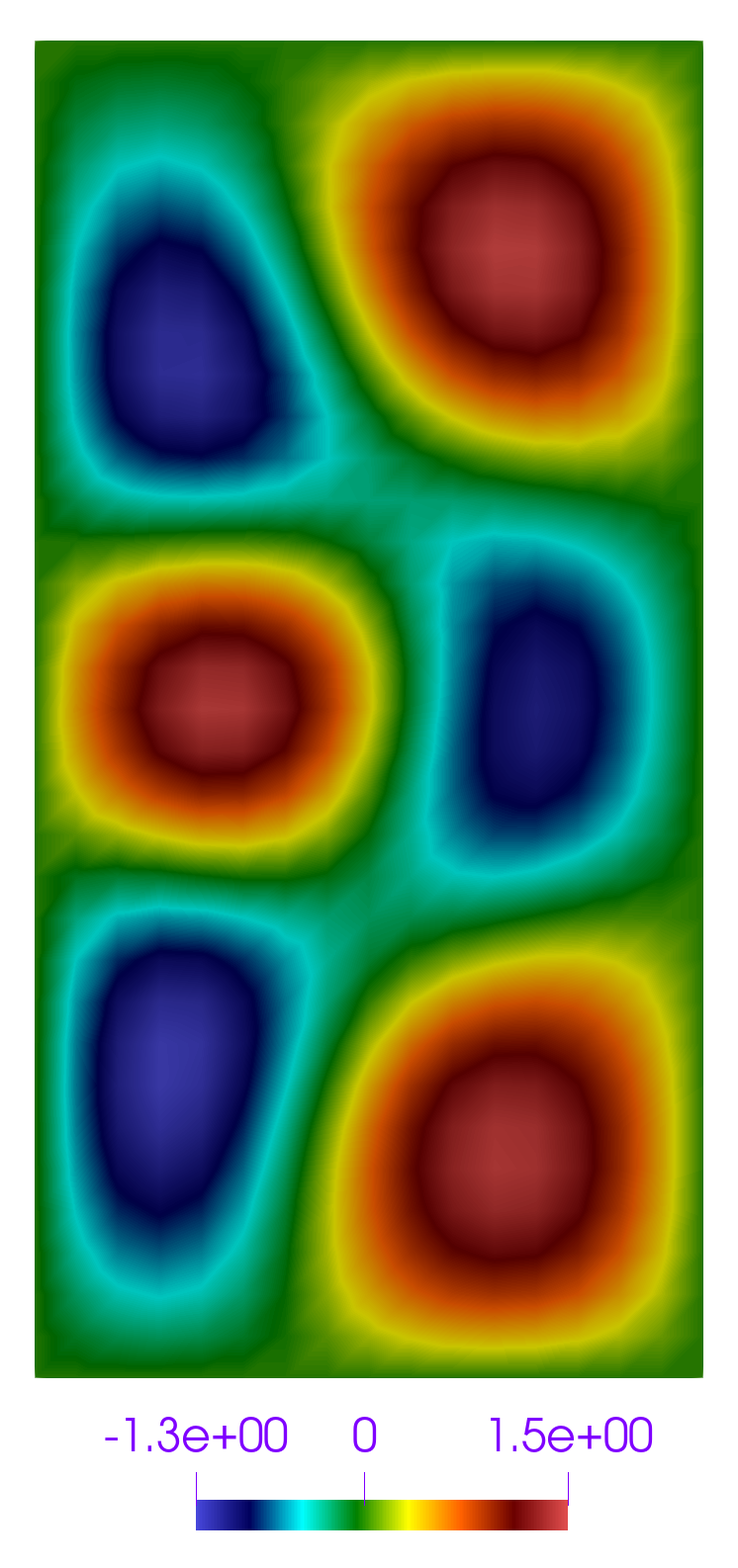}
    \end{overpic}
\put(-48,190){$\xi_7$}
    \begin{overpic}[width=0.2\textwidth]{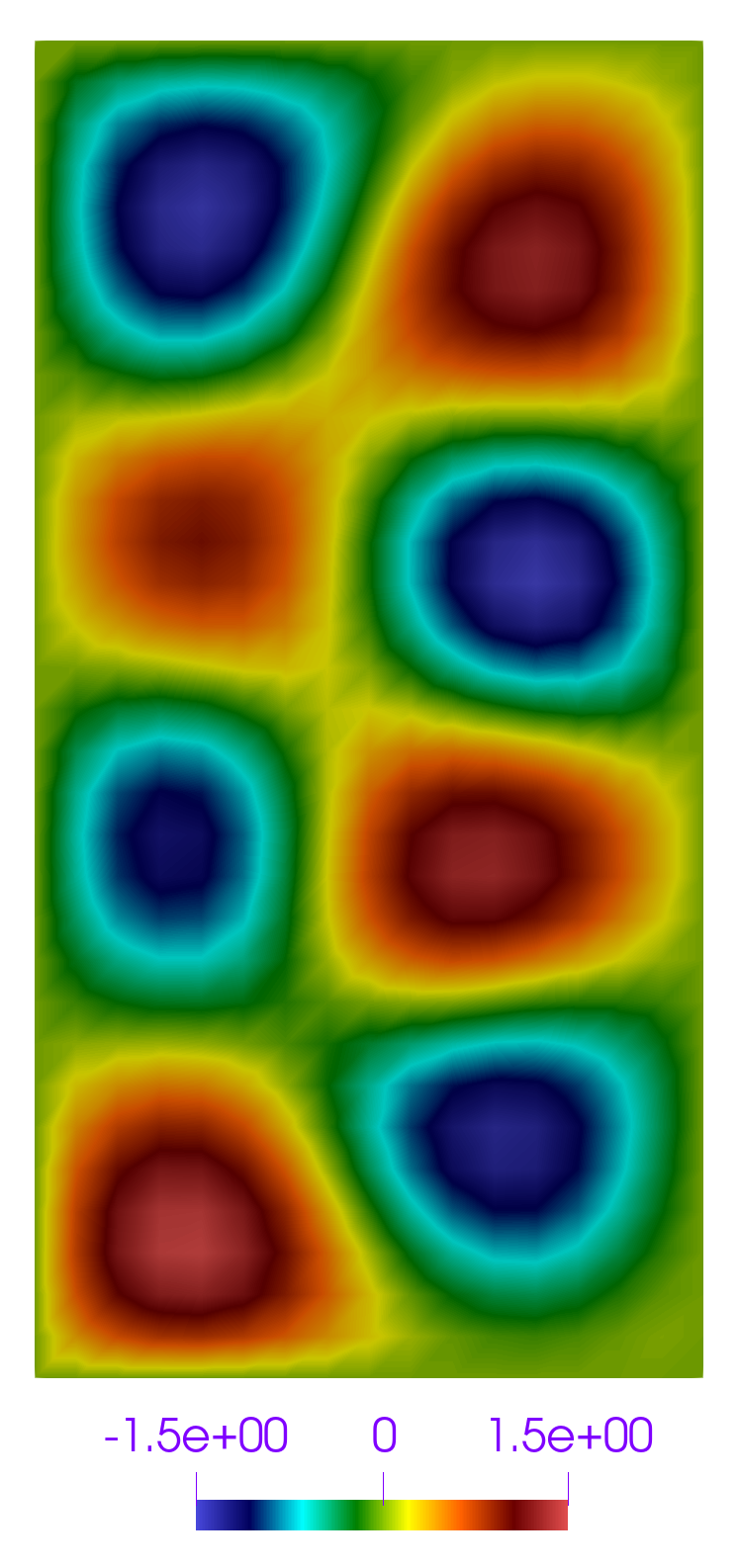}
    \end{overpic}
\put(-48,190){$\xi_{10}$}
    \caption{Case 1: selected POD basis functions for stream function $\psi$.}
    \label{fig:modes}
\end{figure}

For illustration purposes, 
Figure \ref{fig:modes} shows some selected POD basis functions for $\psi$. As expected, the scale of spatial structures becomes smaller and smaller 
as the basis function index increases. 
This is due to the fact that the POD modes are arranged in order of descending energy content. 

\begin{figure}[htb]
    \centering
    \begin{overpic}[width=0.48\textwidth,grid=false]{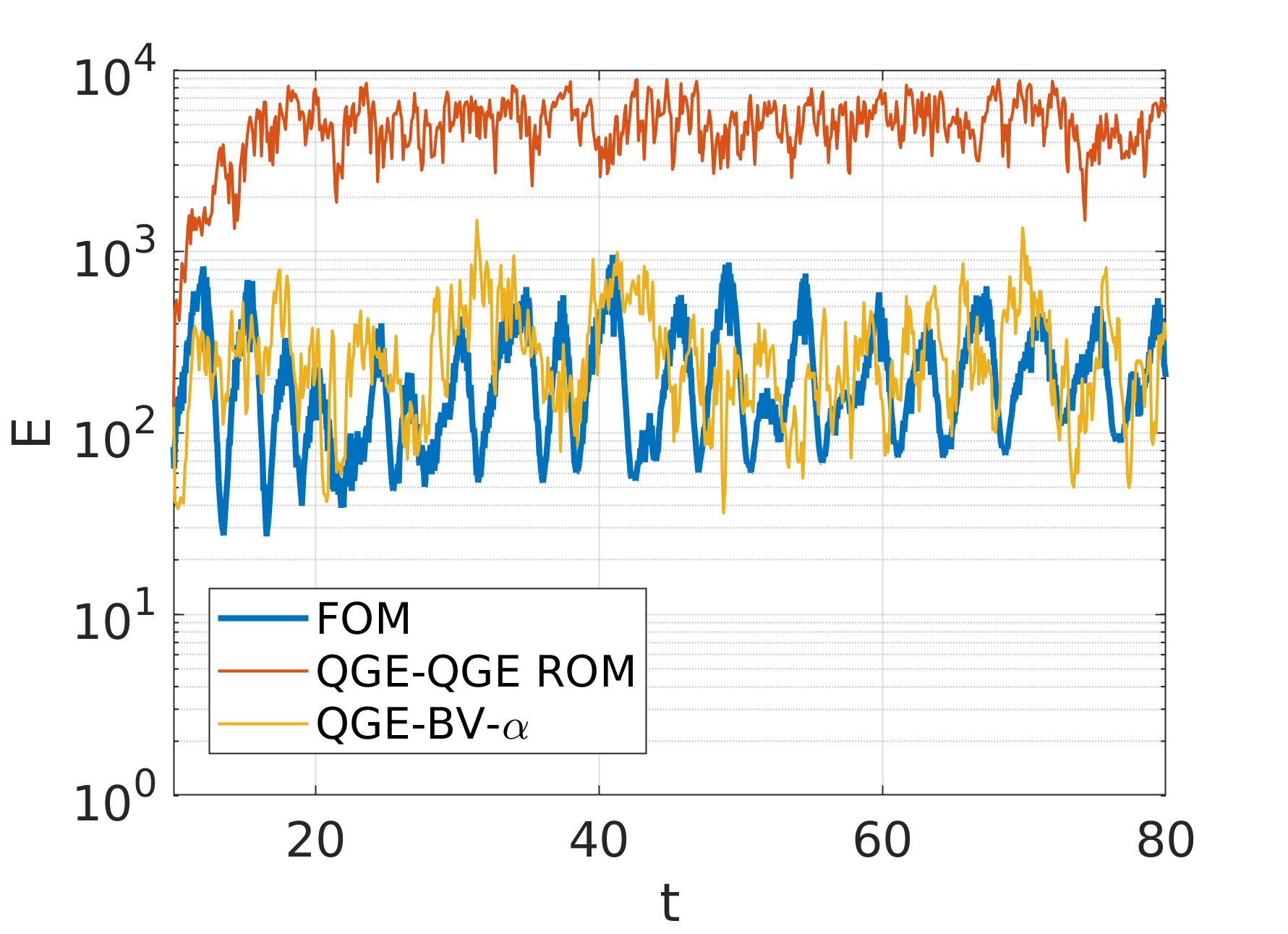}
    \put(42,71){$N_q^r = 10$}
    \end{overpic}
    \begin{overpic}[width=0.48\textwidth]{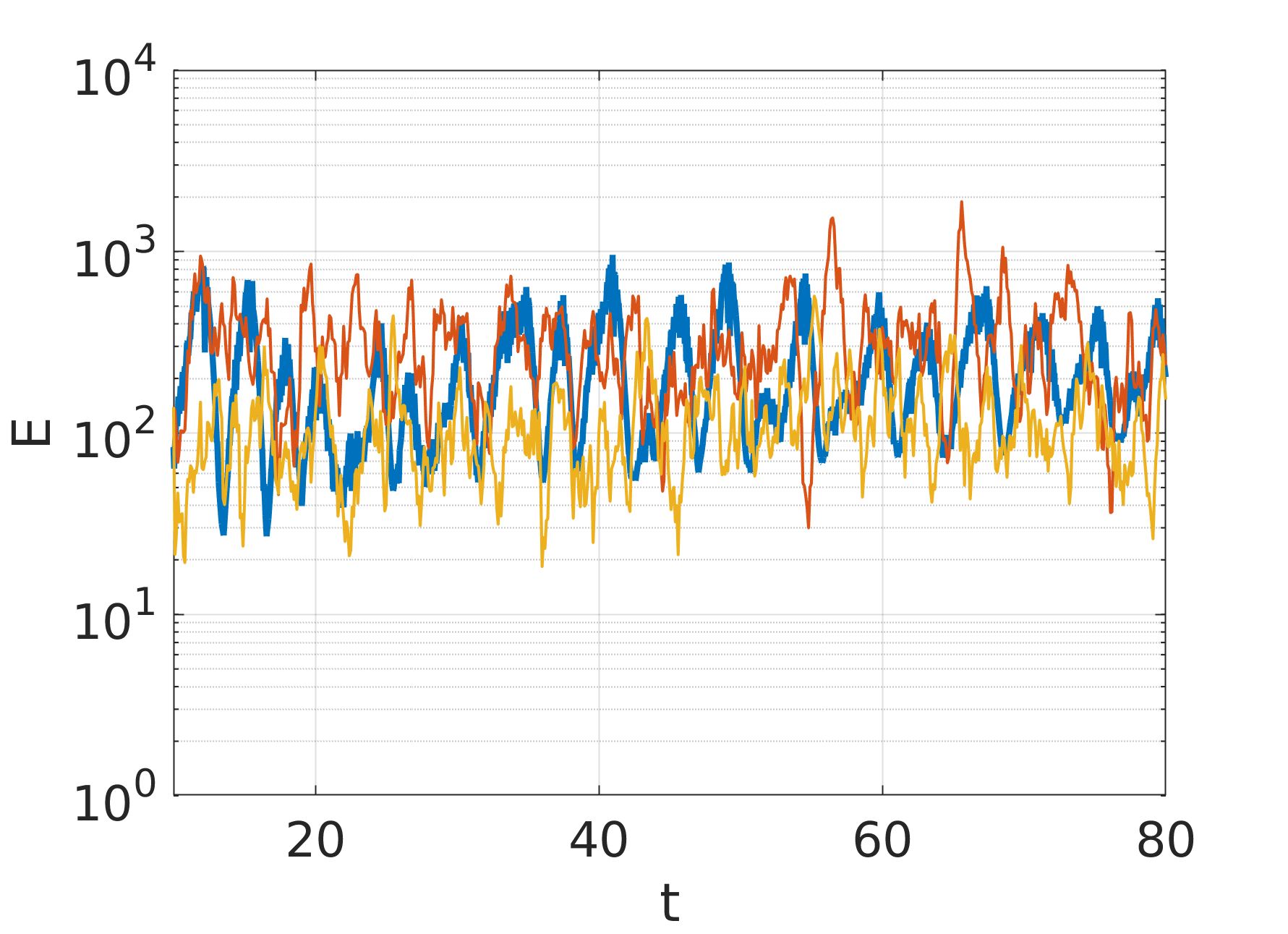}
    \put(42,71){$N_q^r = 20$}
    \end{overpic}
    \begin{overpic}[width=0.48\textwidth]{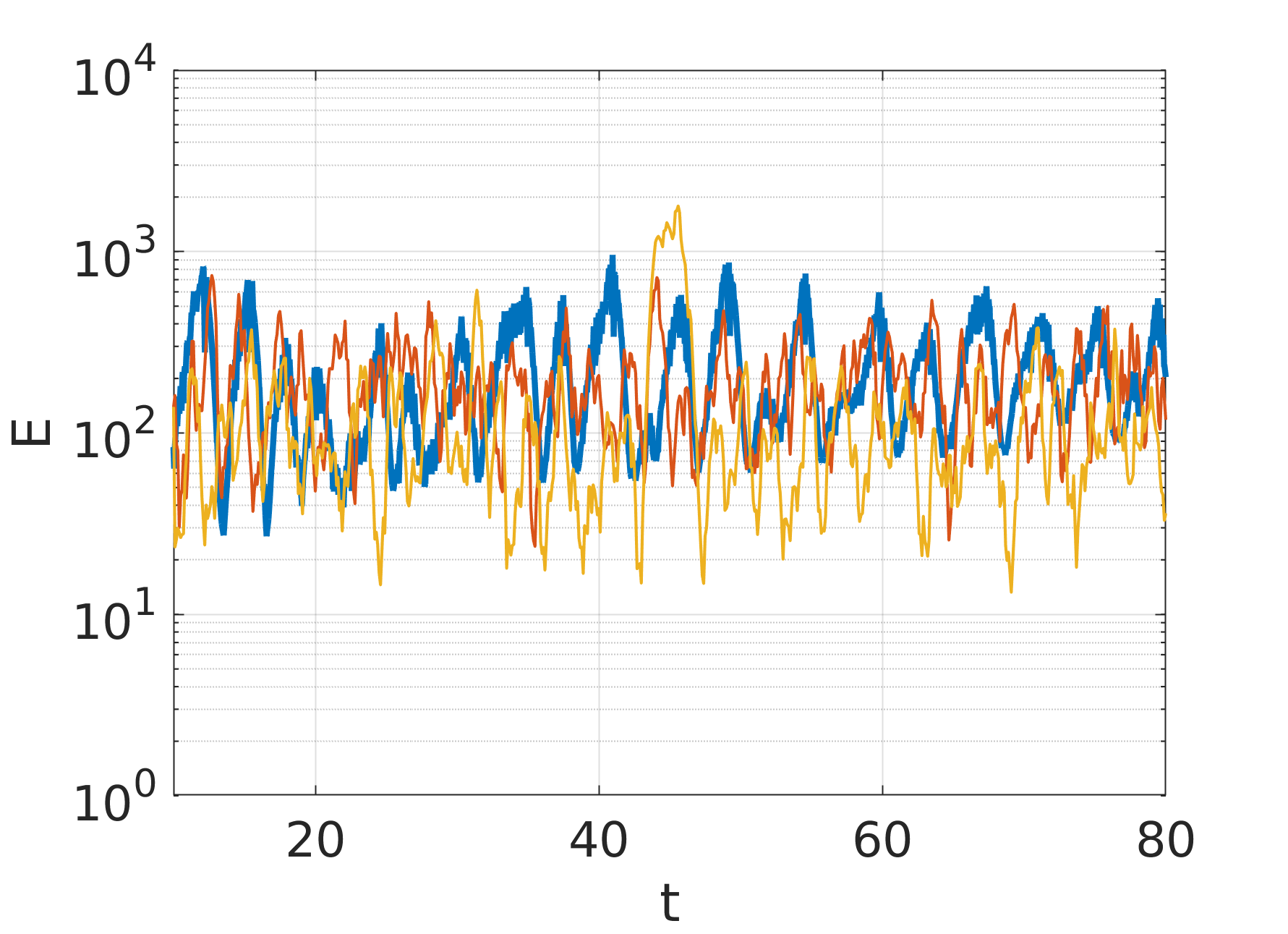}
    \put(42,71){$N_q^r = 30$}
    \end{overpic}
\caption{Case 1: time evolution of the kinetic energy computed by the FOM, QGE-QGE ROM and QGE-BV-$\alpha$ ROM for different numbers of POD basis functions for the potential vorticity: $N_q^r = 10$ (top left), $N_q^r = 20$ (top right) and $N_q^r = 30$ (bottom). $N_\psi^r$ is set to 10. The legend in the top-left panel is common to all the panels.}
    \label{fig:energy_case1}
\end{figure}

Next, we set the energy threshold for the selection
of the stream function eigenvalues to 98\%, which results in 10 modes (i.e., $N^r_\psi = 10$), and let $N_q^r$ vary. In particular, we consider three values of $N_q^r$ that lead to under-resolution at the reduced order level: $N_q^r = 10, 20, 30$ that are needed to retain only
54\%, 65\% and 70\% of the eigenvalue energy, respectively. 
Figure \ref{fig:energy_case1} compares the time evolution of the kinetic energy $E$ \eqref{eq:kin_energy} computed by the FOM and our two ROM approaches for the three values of $N_q^r$.
We observe that the kinetic energy computed by QGE-QGE ROM with $N_q^r = 10$ is much higher than the FOM kinetic energy over the entire time interval of interest. 
Since with $N_q^r = 10$ we only capture 54\% of the eigenvalue energy associated to $q$, such mismatch is to be expected. 
However, if we switch to the QGE-BV-$\alpha$ ROM with the same $N_q^r$, we obtain an average kinetic energy that compares well with average computed by the FOM. By increasing $N_q^r$ to $20$ or $30$, both the QGE-QGE ROM and the QGE-BV-$\alpha$ ROM provide a good prediction of the average kinetic energy.

\begin{figure}[htb]
    \centering
        \begin{overpic}[width=0.19\textwidth,grid=false]{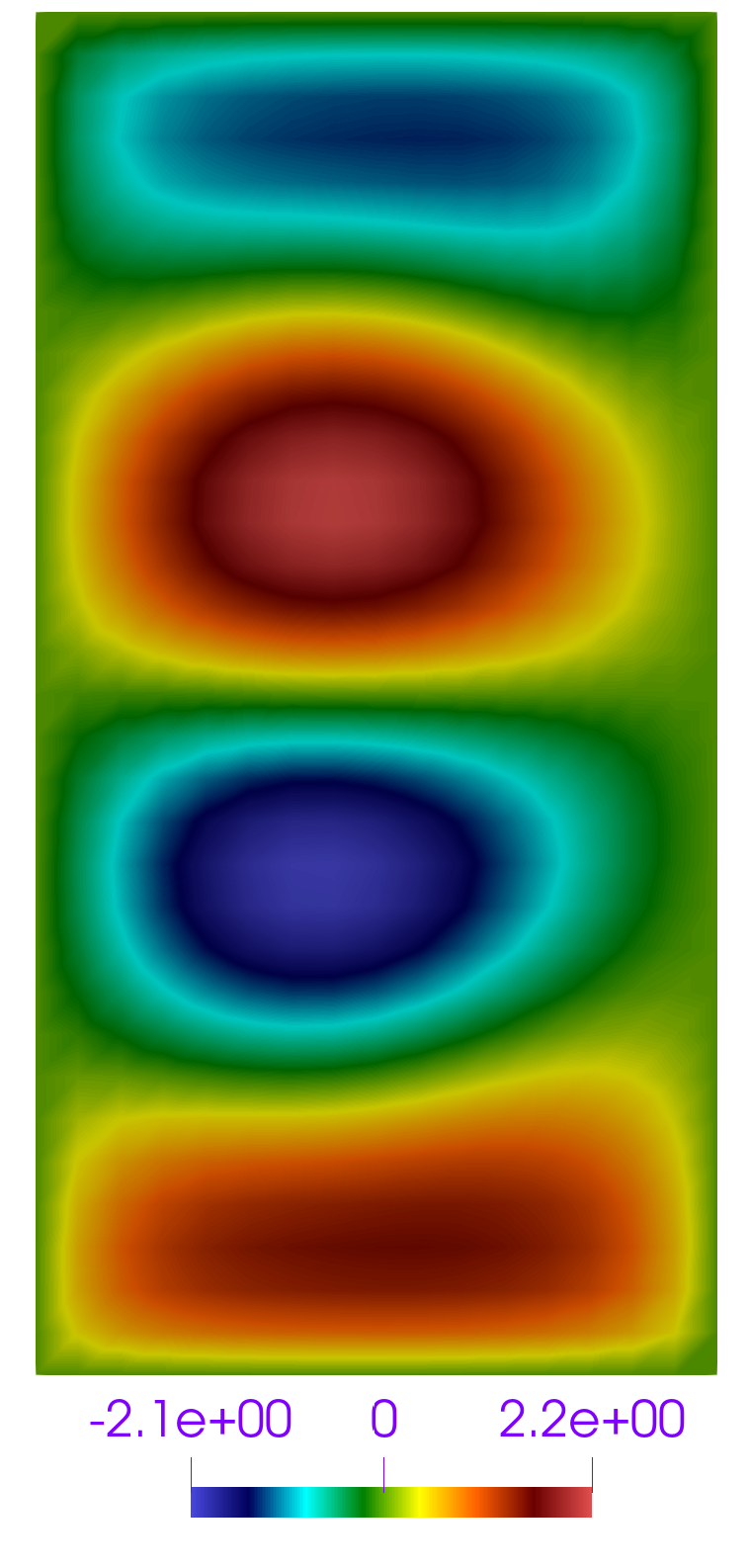}
        \put(17,102){FOM}
        \put(-30,60){QGE-QGE}
        \put(-22,52){ROM}
    \end{overpic}
    \begin{overpic}[width=0.20\textwidth]{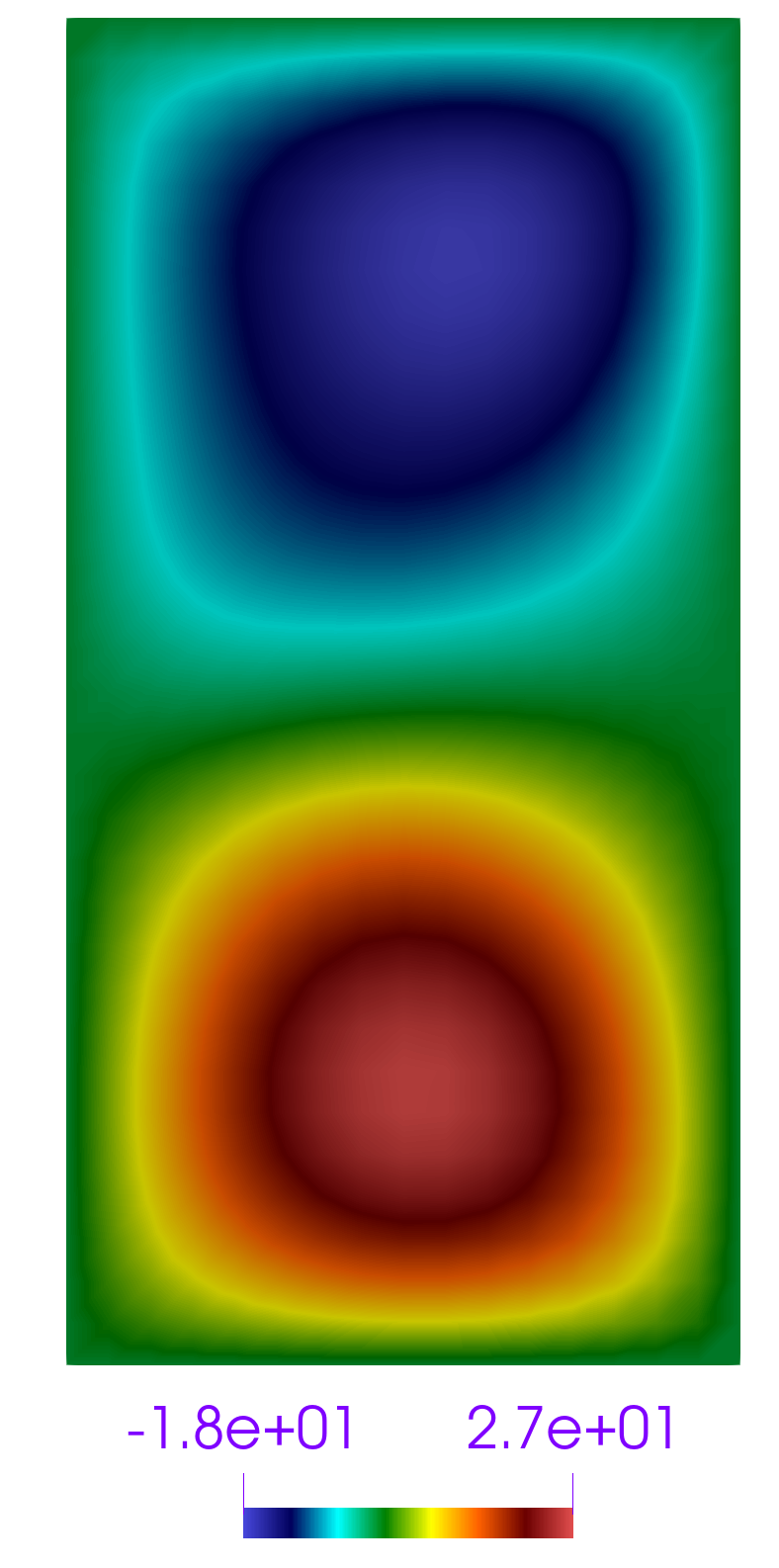}
        \put(17,102){$N_q^r = 10$}
    \end{overpic}
    \begin{overpic}[width=0.2\textwidth]{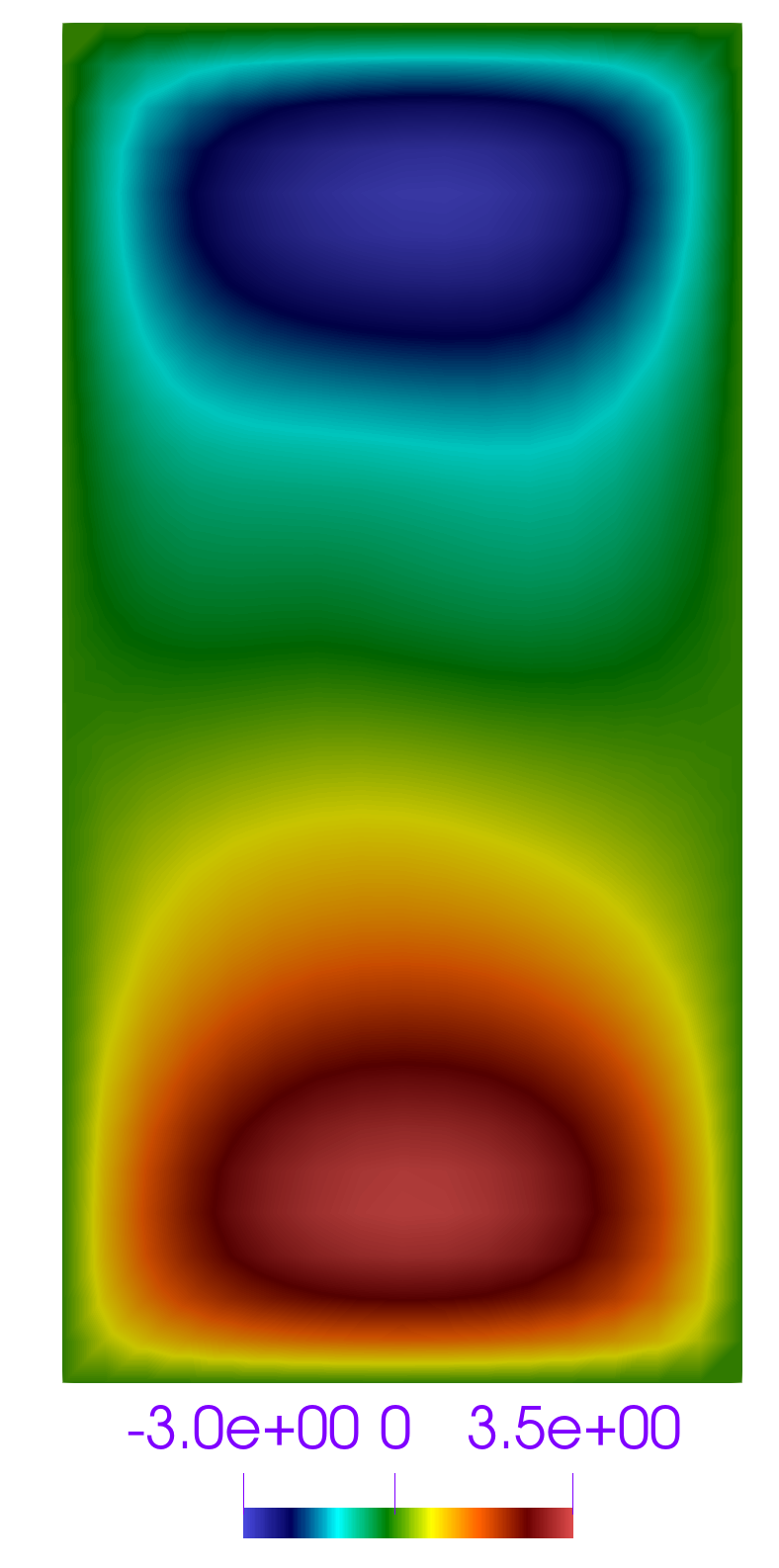}
                \put(17,102){$N_q^r = 20$}
    \end{overpic}
    \begin{overpic}[width=0.2\textwidth]{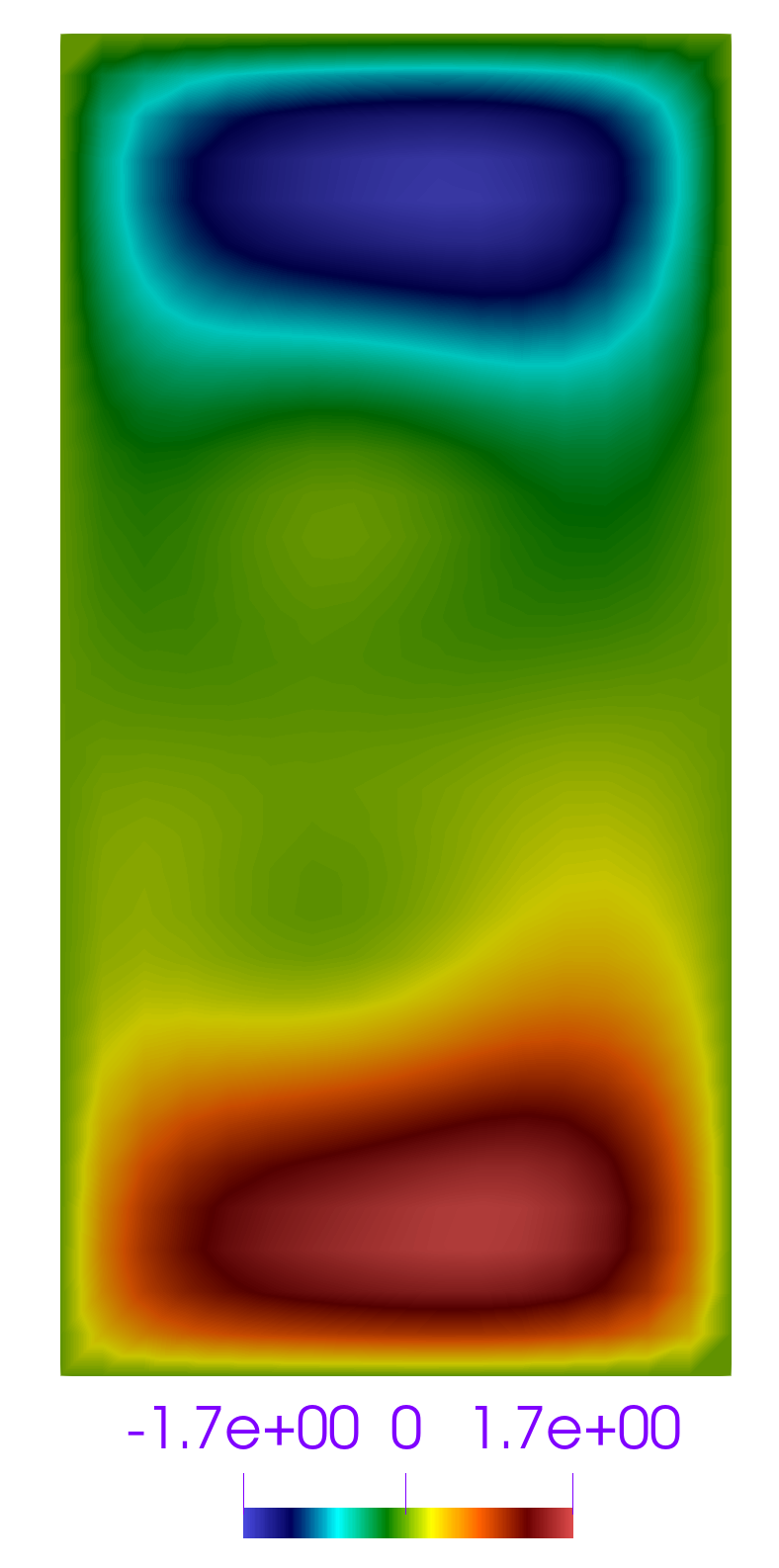}
            \put(17,102){$N_q^r = 30$}
    \end{overpic}
\\
\vskip .5cm
        \begin{overpic}[width=0.19\textwidth]{sections/img/FOM_case1.png}
        \put(17,102){FOM}
        \put(-30,60){QGE-BV-$\alpha$}
        \put(-22,52){ROM}
    \end{overpic}
    \begin{overpic}[width=0.2\textwidth]{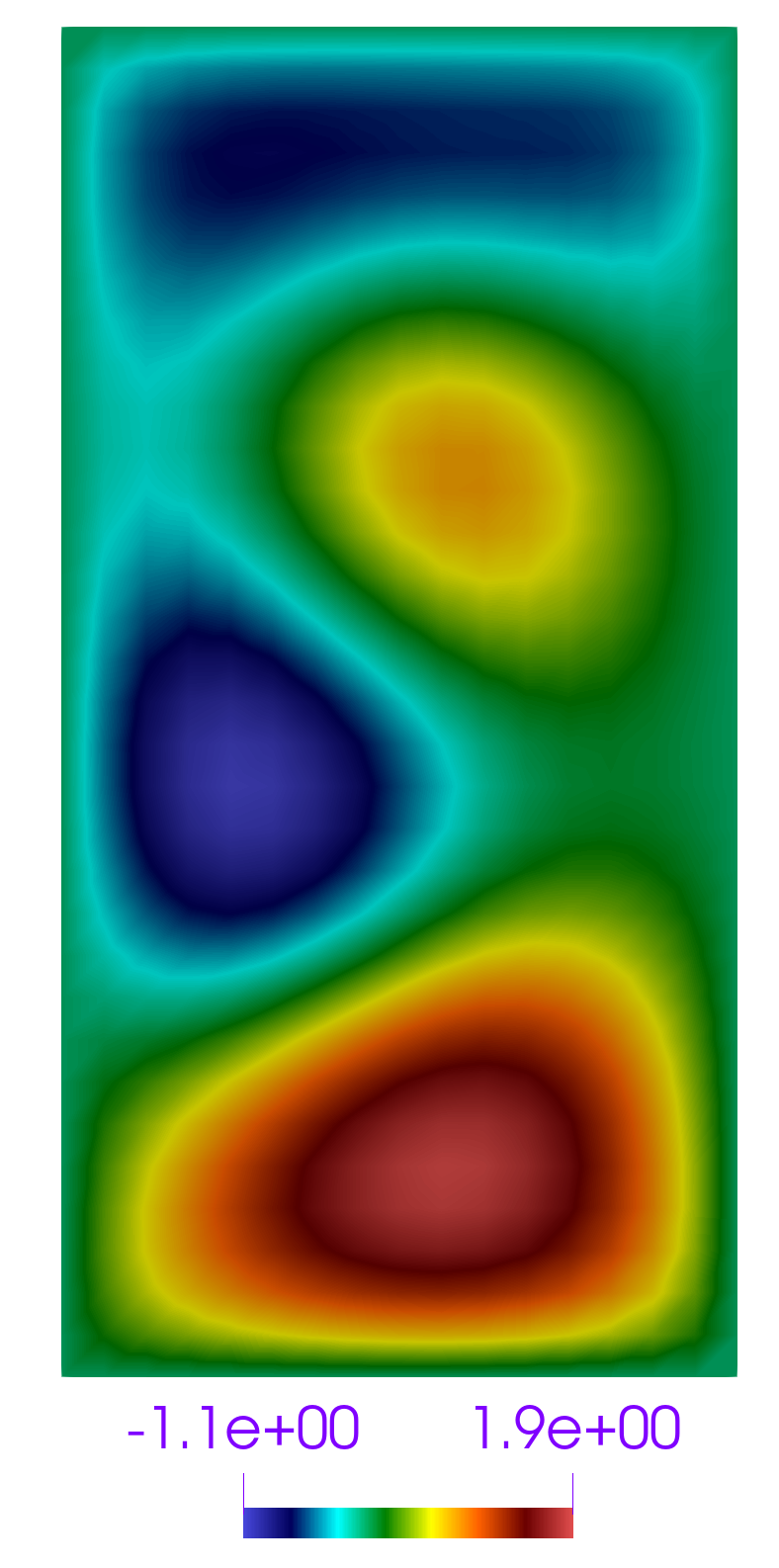}
            \put(17,102){$N_q^r = 10$}
    \end{overpic}
    \begin{overpic}[width=0.19\textwidth]{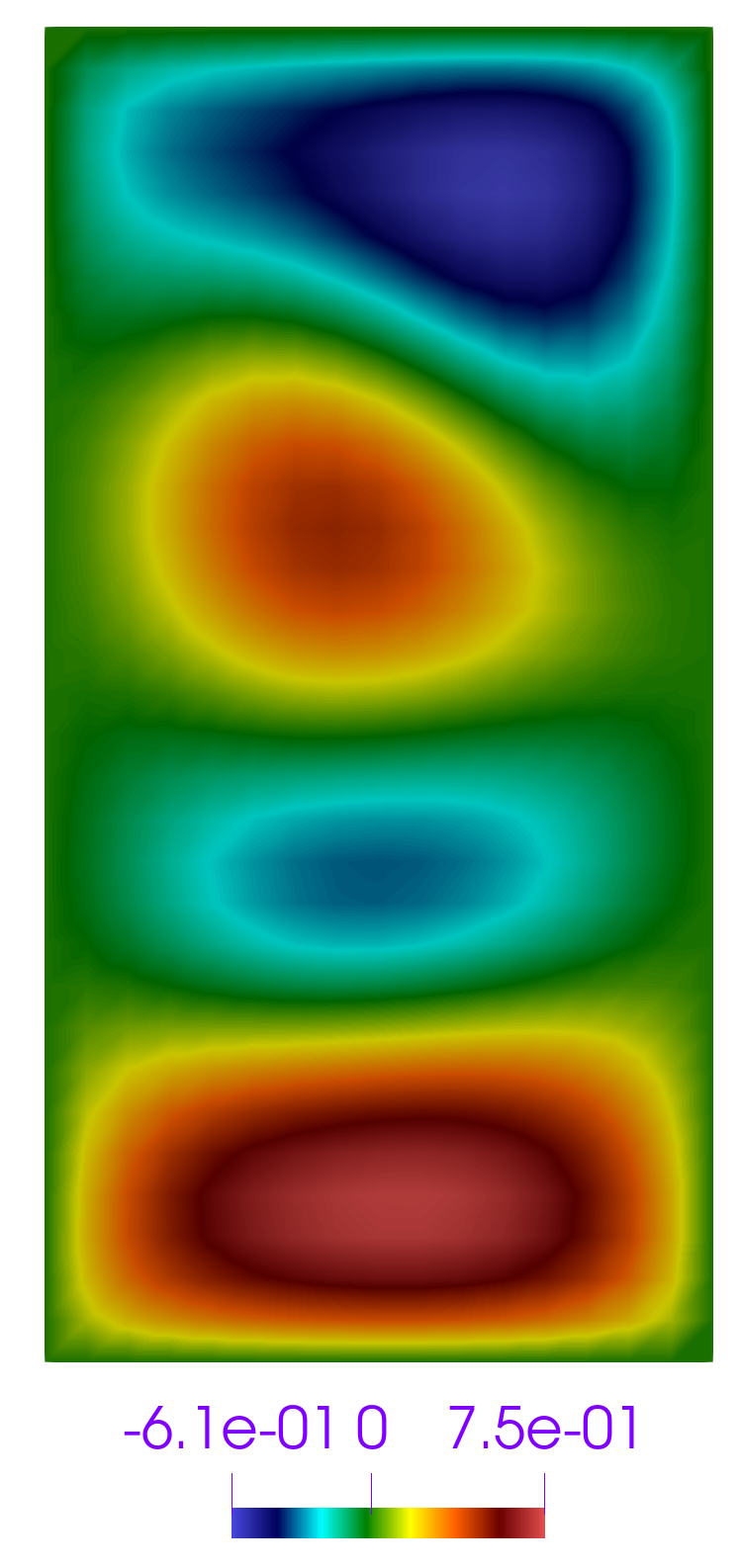}
                \put(17,102){$N_q^r = 20$}
    \end{overpic}
    \begin{overpic}[width=0.19\textwidth]{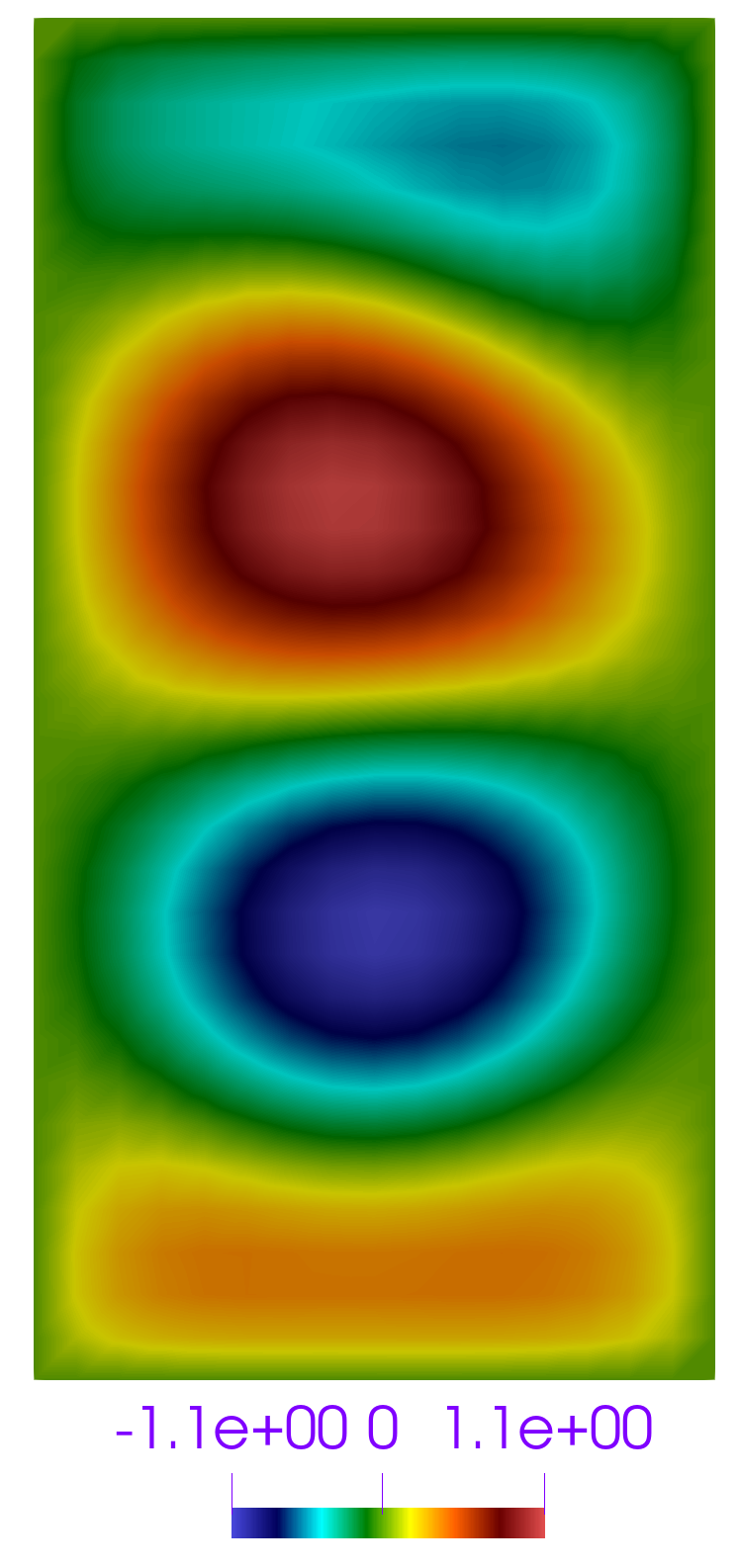}
                \put(17,102){$N_q^r = 30$}
    \end{overpic}
    \caption{Case 1: comparison of $\widetilde{\psi}$ computed by the FOM and the QGE-QGE ROM (top row)
    or the QGE-BV-$\alpha$ ROM (bottom row) for different numbers of POD vorticity modes $N_q^r$. $N_\psi^r$ is set to 10.}
    \label{fig:psiROMstand}
\end{figure}
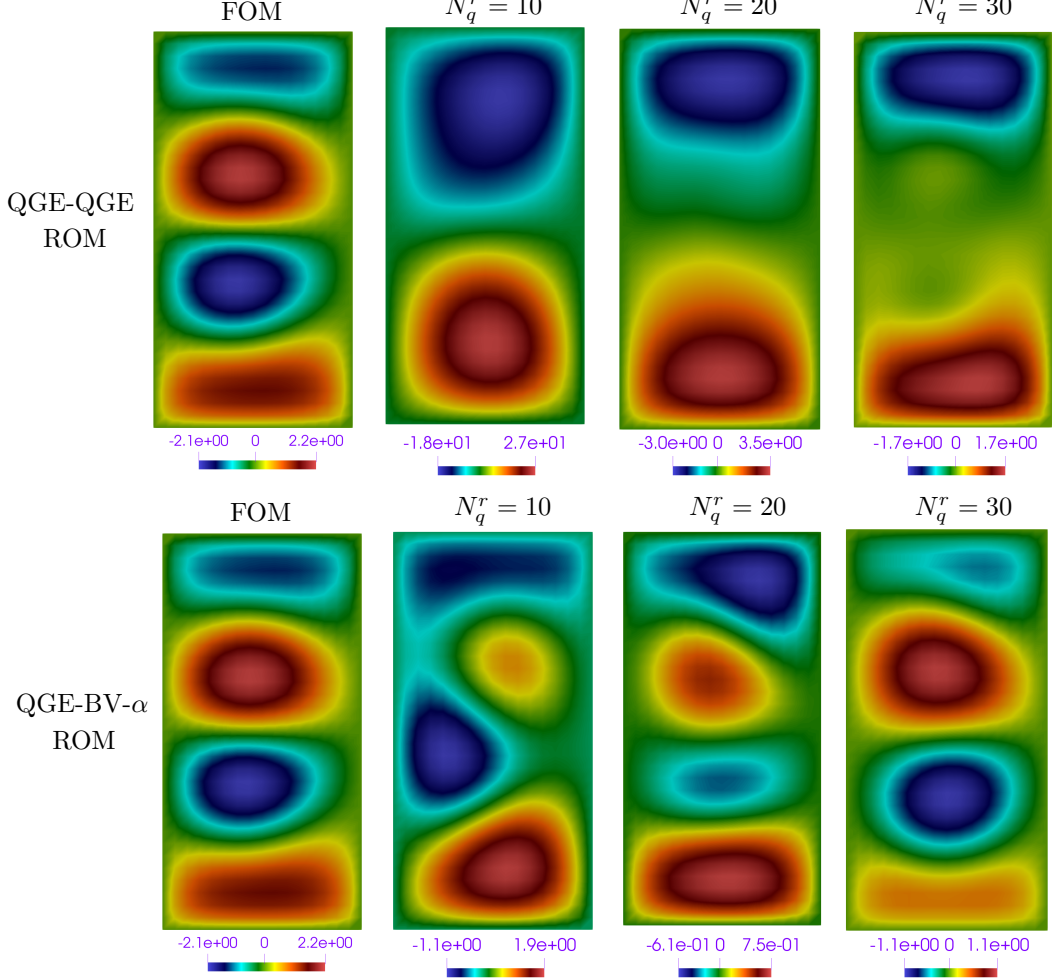

For further comparison, Figure \ref{fig:psiROMstand}
reports the time-averaged stream function \eqref{eq:time_psi} computed by the FOM and QGE-QGE and QGE-BV-$\alpha$ ROMs for the same values of $N_q^r$ as in Fig.~\ref{fig:energy_case1}. As expected from Fig.~\ref{fig:energy_case1} (top-left panel), $\widetilde{\psi}$ computed by the QGE-QGE ROM with 
$N_q^r = 10$ is highly inaccurate. See the second panel in the top row of Fig.~\ref{fig:psiROMstand}, which shows only two gyres instead of four. However, even when the $N_q^r$ is increased to 20 and 30 the QGE-QGE ROM fails to reproduce the four-gyre pattern
despite the fact that the average kinetic energy is well captured.
It is interesting to note that these reduced order solution computed with $N_q^r = 10, 20, 30$ look very similar to full order solutions computed
with the QGE model on a severely under-refined mesh. See Fig.~1 in \cite{Girfoglio1}, which were obtained with a $4 \times 8$ mesh. This is evidence of the analogy between an under-resolved ROM and an under-resolved FOM.
The second panel on the bottom row of Fig.~\ref{fig:psiROMstand} shows that the
QGE-BV-$\alpha$ ROM is able to recover the four-gyre pattern already with $N_q^r = 10$, although obviously the solution is not accurate. 
As $N_q^r$ is increased, QGE-BV-$\alpha$ ROM provides solutions that get closer and closer to the FOM solution. Even with $N_q^r = 30$ though, 
the magnitude of $\widetilde{\psi}$ computed by the QGE-BV-$\alpha$ ROM is smaller than it should be. See the last panel on the bottom row of Fig.~\ref{fig:psiROMstand}. We suspect that this is due to the use of a linear filter. In fact, when a linear filter is adopted at the full order level, the solutions are characterized by over-diffusion since the filter is not selective. It is reasonable to expect a similar behavior when a linear filter is used at the reduced order level. 

To make the comparison between QGE-QGE and QGE-BV-$\alpha$ ROMs more quantitative, we report in Table \ref{tab:errors} the $L^2$ errors \eqref{eq:normerror}. We see 
for any value of $N_q^r$ the $L^2$ error obtained with the QGE-BV-$\alpha$ ROM is smaller than the
$L^2$-norm error obtained with the QGE-BV-$\alpha$ ROM.
This is particularly evident when $N_q^r = 10$:
the QGE-BV-$\alpha$ ROM provides an error about 15 times smaller than the QGE-QGE ROM.

\begin{table}[htb]
\centering
\begin{tabular}{|c|c|c|c|}
\multicolumn{4}{c}{} \\
\cline{1-4}
$N_q^r$ & \% of energy content & $\varepsilon_{QGE-QGE}$ & $\varepsilon_{QGE-BV-\alpha}$\\
\hline
10 & 54\% & 1.2e+01 & 8.1e-01 \\
20 & 65\% & 1.7e+00 & 7.7e-01\\
30 & 70\% & 9.2e-01 & 6.1e-01\\
\hline
\end{tabular}
\caption{Case 1: $L^2$ error \eqref{eq:normerror} given by QGE-QGE-ROM and QGE-BV-$\alpha$ ROM for $N_q^r = 10, 20, 30$ and $N_\psi = 10$.}
\label{tab:errors}
\end{table}

\begin{figure}
        \begin{overpic}[width=0.182\textwidth]{sections/img/FOM_case1.png}
    \end{overpic}
    \put(-58,180){FOM}
        \begin{overpic}[width=0.19\textwidth]{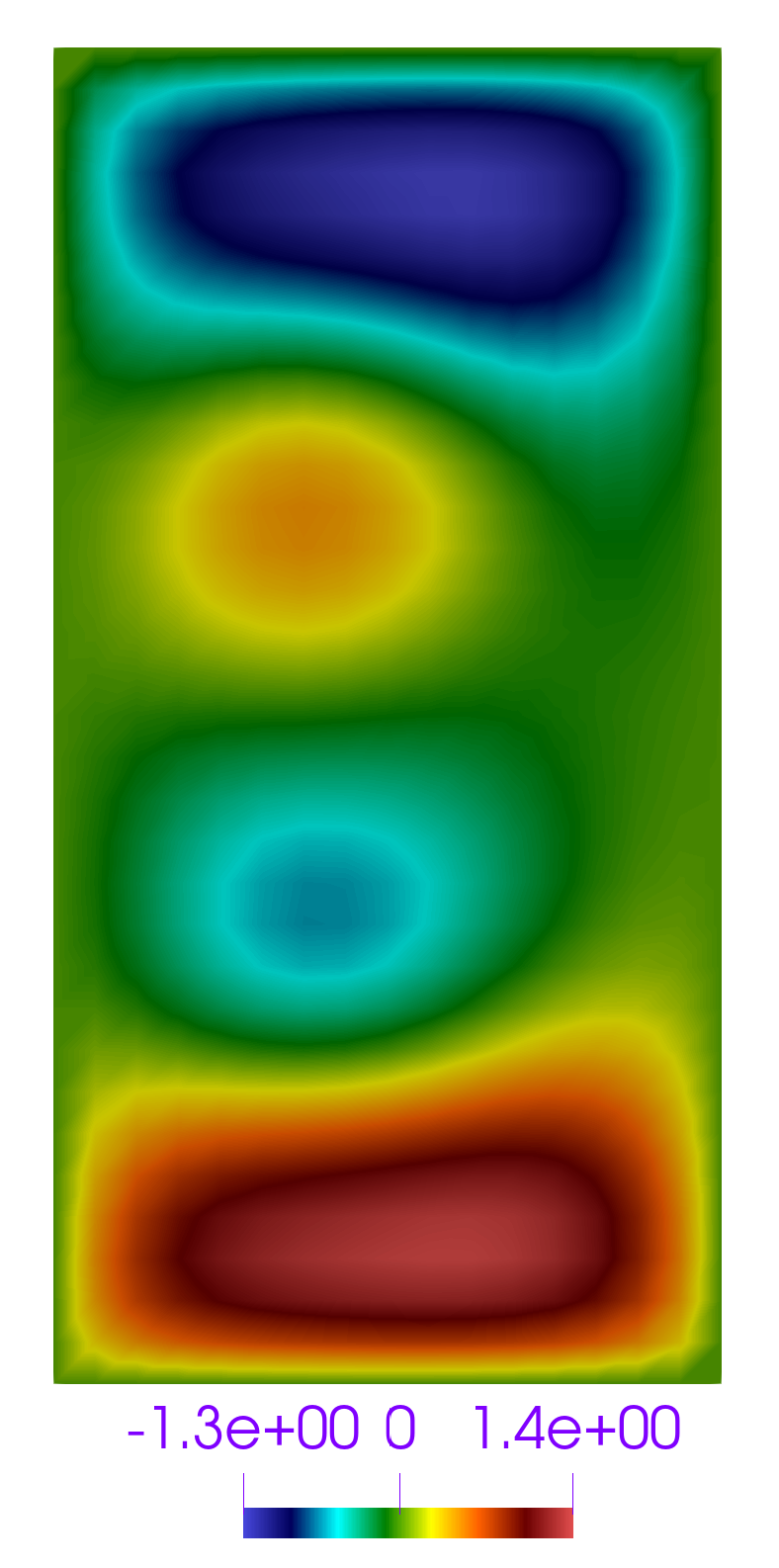}
    \end{overpic}
\put(-58,180){$N_q^r = 40$}
    \caption{Case 1: $\widetilde{\psi}$ computed by the FOM and the QGE-QGE ROM with $N_q^r = 40$ and $N_\psi^r = 10$.}
    \label{fig:psiROMalpha}
\end{figure}

If one wanted to retrieve the four-gyre pattern
with the QGE-QGE ROM, $N_q^r$ has to be increased to $40$, which corresponds to retaining 74\% of the eigenvalue energy. See Fig.~\ref{fig:psiROMalpha}. This is in agreement with \cite{QGE_review}. However, despite displaying the correct pattern, the  $\widetilde{\psi}$ computed by the QGE-QGE ROM is still far from being accurate.

Finally, we provide a comment on the efficiency of our ROM approaches. 
Table \ref{tab:time} reports the CPU time required by the QGE-QGE ROM and the QGE-BV-$\alpha$ ROM with $N_q^r = 10$ and $N_\psi = 10$ and the relative speed-ups with respect to the CPU time a FOM simulation (506 s). We note that by CPU time for the ROMs we mean the total time needed to solve the linear systems that yield the modal coefficients, i.e. the CPU time taken by the online phase only. 
A first observation on Table \ref{tab:time} is that the higher accuracy of the QGE-BV-$\alpha$ ROM comes with an increased of about 50\% in computational time with respect to the QGE-QGE ROM. Furthermore, we observe that, despite the overall low accuracy of the ROM solutions for $N_q^r = 10$ and $N_\psi = 10$, the speed-up is not particularly encouraging. If $N_q^r$ is increased to 30 to gain accuracy, the speed-up deteriorates further. One reason we have identified for this poor reduction of the computational cost is the choice of a coarse mesh (i.e., 16 $\times$ 32). If we were to use a much finer mesh (e.g., 256 $\times$ 512 as in \cite{Girfoglio1}), we would see more important computational savings. However, that would entail a much more onerous offline phase. Indeed, a FOM simulation with mesh 256 $\times$ 512 takes about 10 hours, leading to relative speed-ups two order of magnitude larger than the ones in Tables \ref{tab:time}.

\begin{table}[htb]
\centering
\begin{tabular}{|c|c|c|}
\hline
& $\varepsilon_{QGE-QGE}$ & $\varepsilon_{QGE-BV-\alpha}$\\
\hline
CPU time & 105 s & 165 s \\
\hline
speed-up & 4.8 & 3.1\\
\hline
\end{tabular}
\caption{Case 1: CPU time required by 
QGE-QGE-ROM and QGE-BV-$\alpha$ ROM with $N_q^r = 10$ and $N_\psi = 10$ and relative speed-up with respect to the CPU time required by the FOM simulation (506 s).}
\label{tab:time}
\end{table}

\subsection{Results for Case 2}


The obvious effect of dealing with a smaller Kolmogorov scale (i.e., smaller $Re$) is the need for a finer mesh. Indeed, as shown in \cite{Girfoglio1} a simulation obtained with the QGE model and mesh 16 $\times$ 32 does not  provide a physical solution for Case 2. By increasing the resolution to 32 $\times$ 64, we do obtain a physical solution in terms of patterns and magnitudes despite the fact that a DNS would require an even finer mesh. Time step is set to $\Delta t = 1e-4$, as in Case 1.

\begin{figure}[htb]
    \centering
    \begin{overpic}[width=0.5\textwidth]{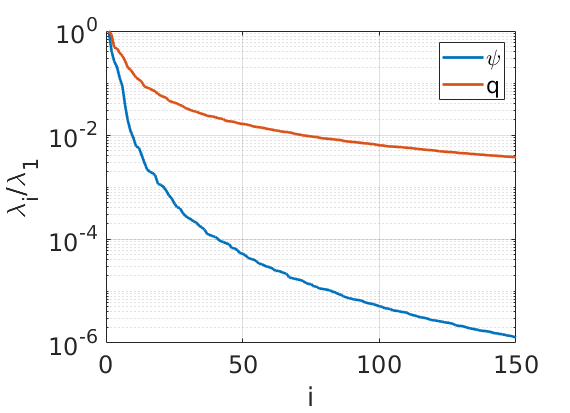}
    \end{overpic}
    \caption{Case 2:
    eigenvalue decay for the stream
function and the vorticity.}
    \label{fig:eig2}
\end{figure}

The eigenvalue decay for the stream function and the  potential vorticity for Case 2 is shown in
Figure \ref{fig:eig2}.  Like for Case 1, the eigenvalue decay is very different for the two variables and it is much slower
for $q$. However, there is one important difference: the eigenvalue decay is slightly faster for both variables in Case 2 than in Case 1. For example,
for $N_q = 30$ we are able to capture 76\% of the cumulative eigenvalue energy for $q$, instead of 70\%. As we will see, this leads to better ROM reconstructions. Thus, while Case 2 seems more challenging than Case 1 at the FOM level, it is less so at the ROM level.

\begin{figure}[htb]
    \centering
    \begin{overpic}[width=0.5\textwidth,grid=false]{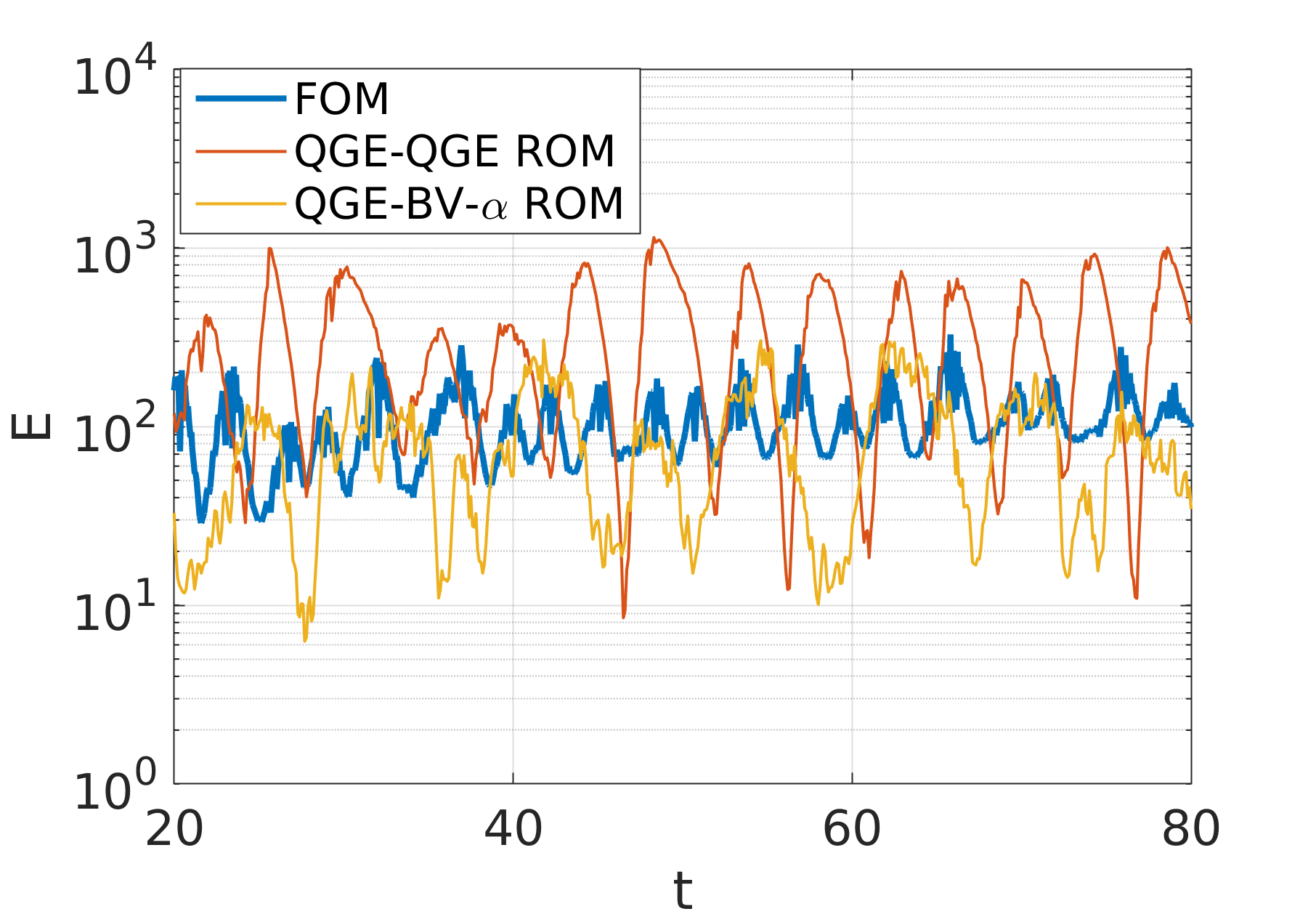}
    \put(42,69){$N_q^r = 10$}
    \end{overpic}
    \begin{overpic}[width=0.48\textwidth]{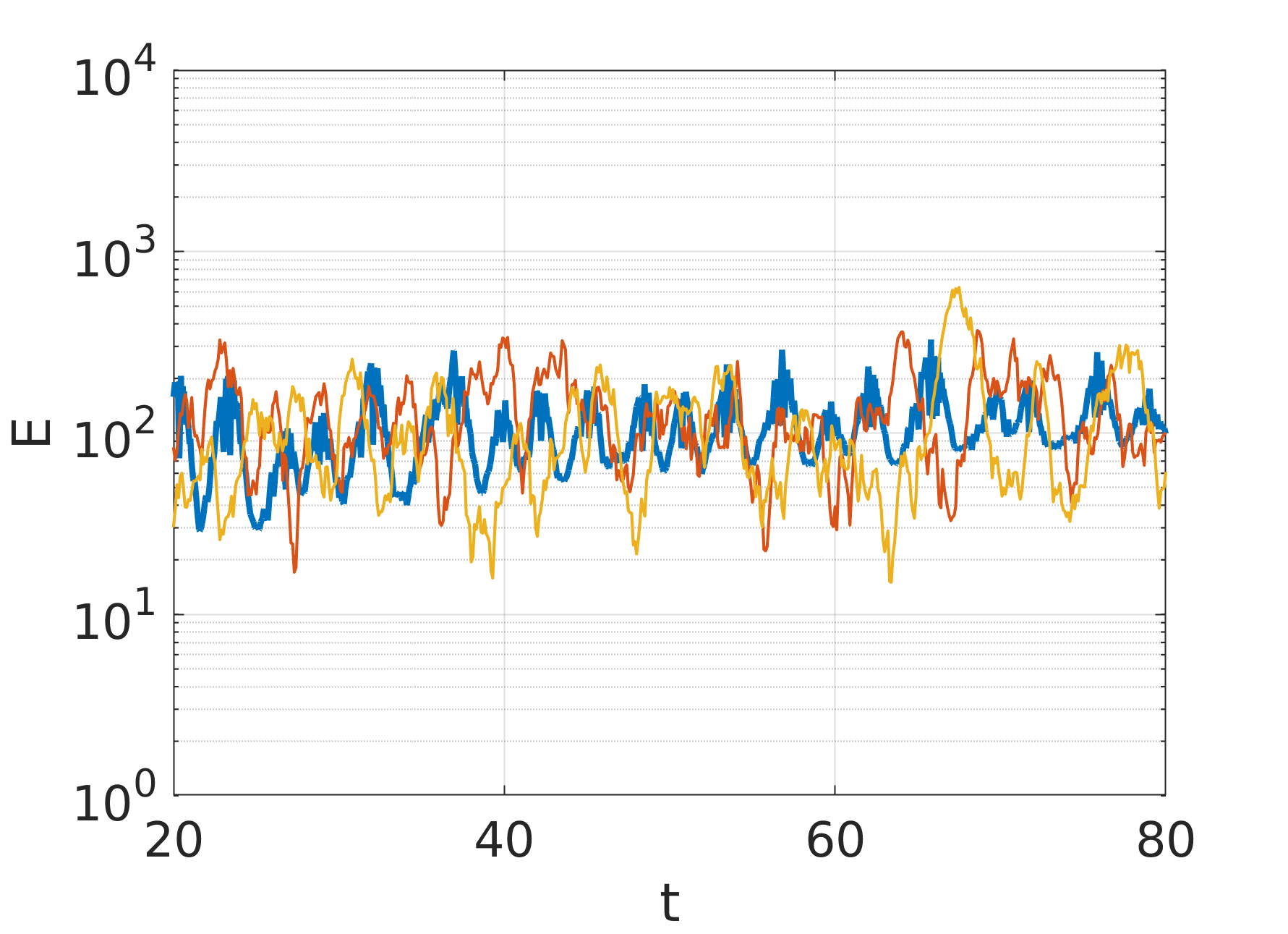}
    \put(42,71){$N_q^r = 20$}
    \end{overpic}
    \begin{overpic}[width=0.48\textwidth]{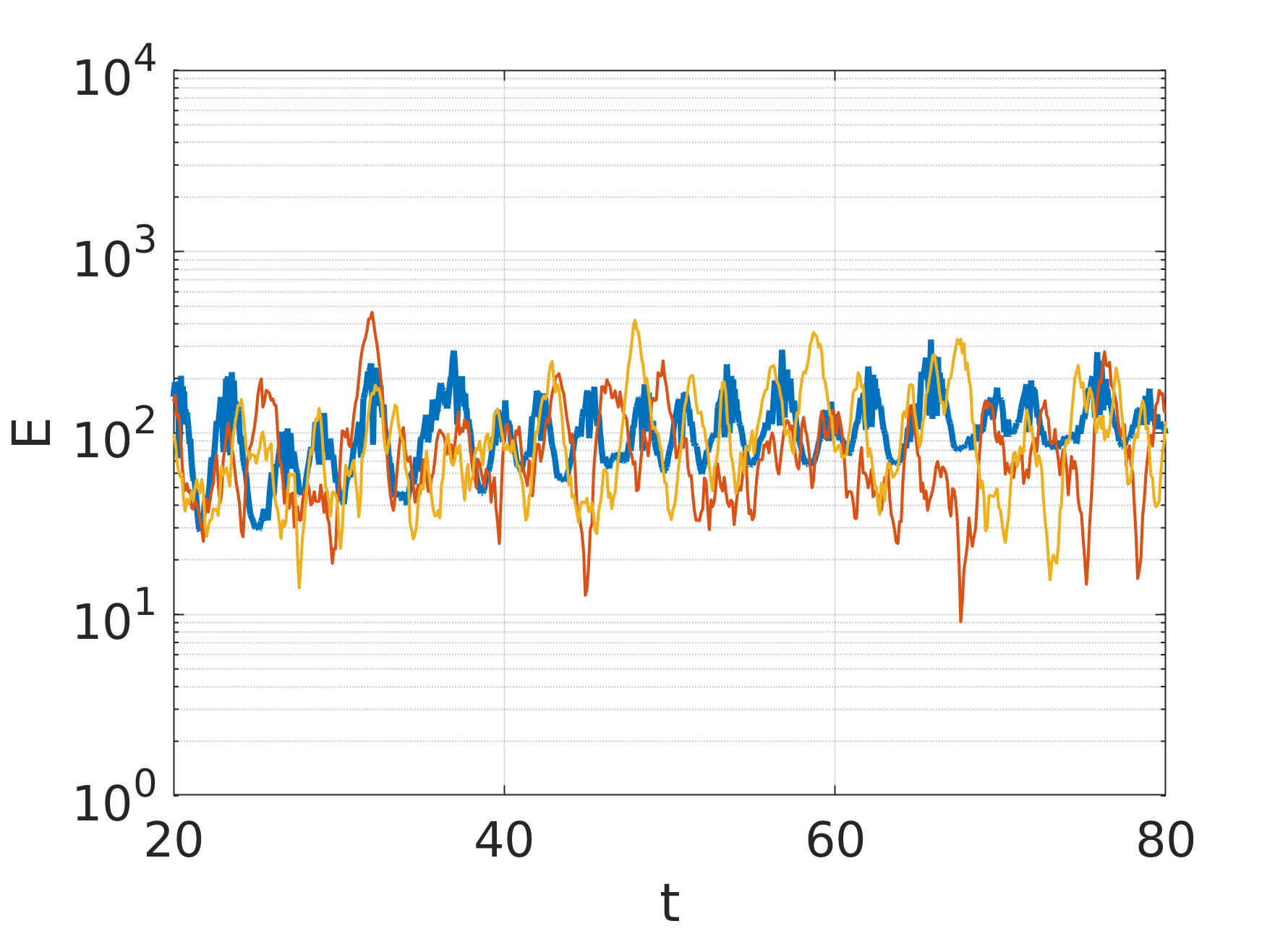}
    \put(42,71){$N_q^r = 30$}
    \end{overpic}
\caption{Case 2: time evolution of the kinetic energy computed by the FOM, QGE-QGE ROM and QGE-BV-$\alpha$ ROM for different numbers of POD basis functions for the potential vorticity: $N_q^r = 10$ (top left), $N_q^r = 20$ (top right) and $N_q^r = 30$ (bottom). $N_\psi^r$ is set to 10. The legend in the top-left panel is common to all the panels.}
    \label{fig:energy_case2}
\end{figure}

\begin{figure}[htb]
    \centering
    \begin{overpic}[width=0.192\textwidth]{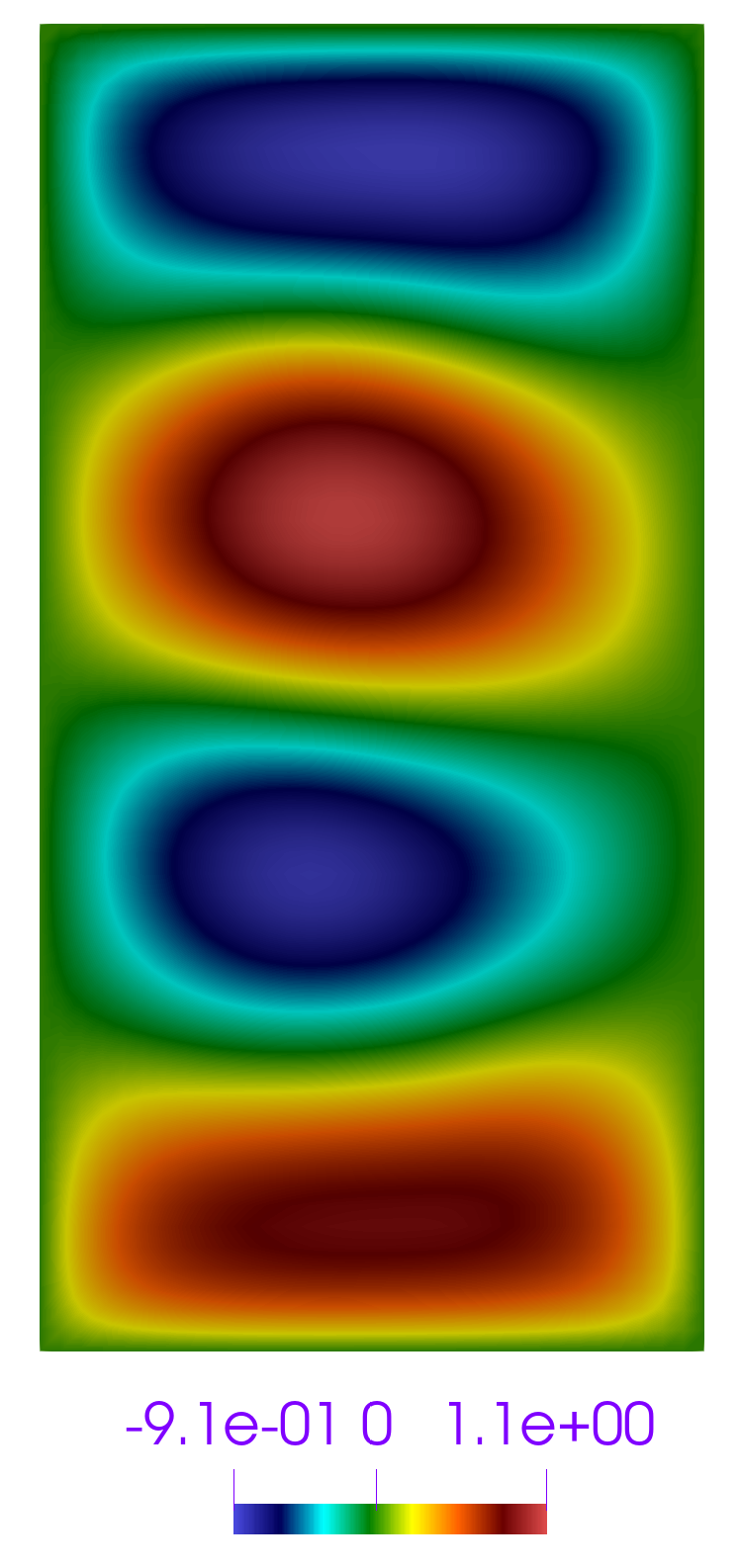}
    \end{overpic}
    \put(-58,180){FOM}
    \begin{overpic}[width=0.192\textwidth]{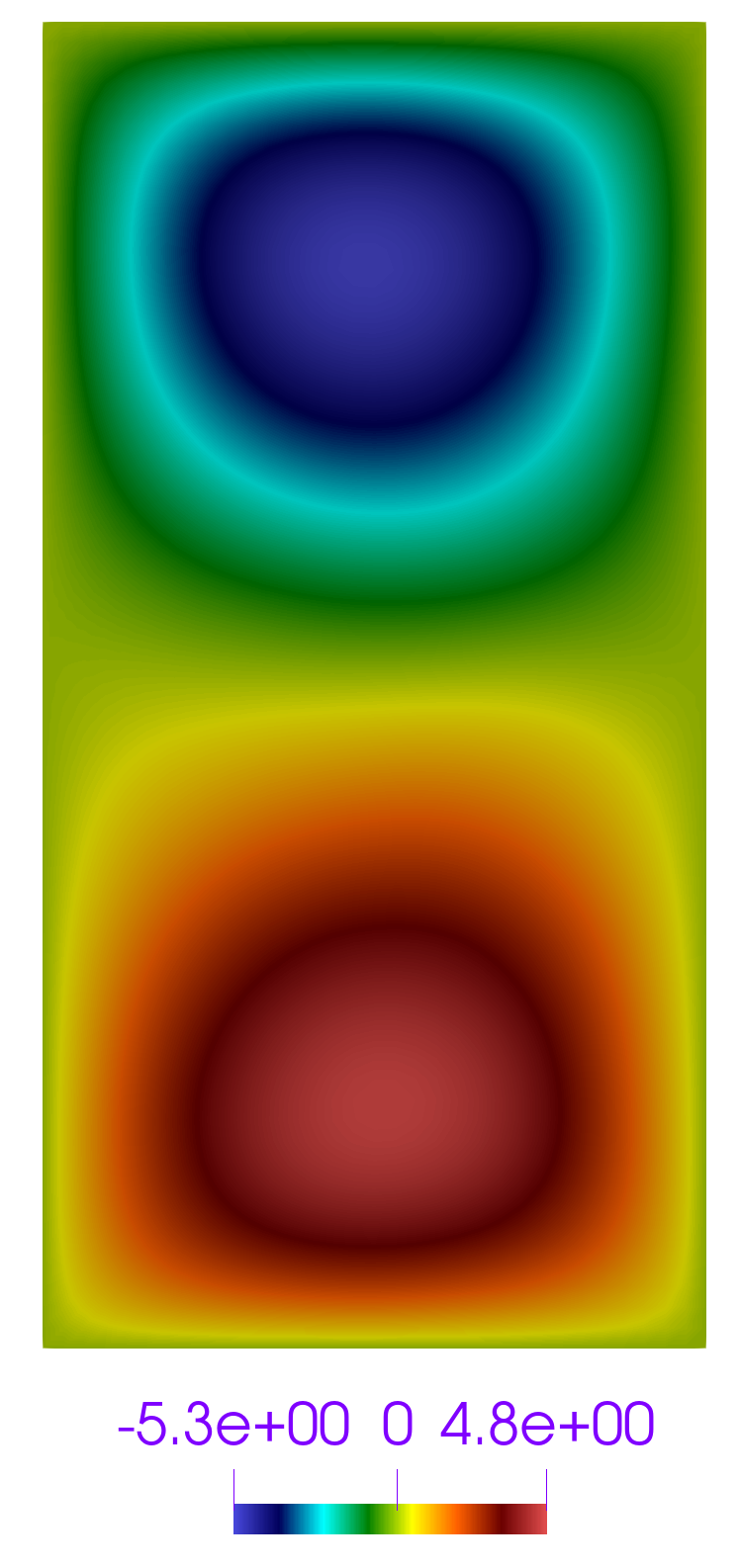}
    \end{overpic}
\put(-58,180){$N_q^r = 10$}
    \begin{overpic}[width=0.192\textwidth]{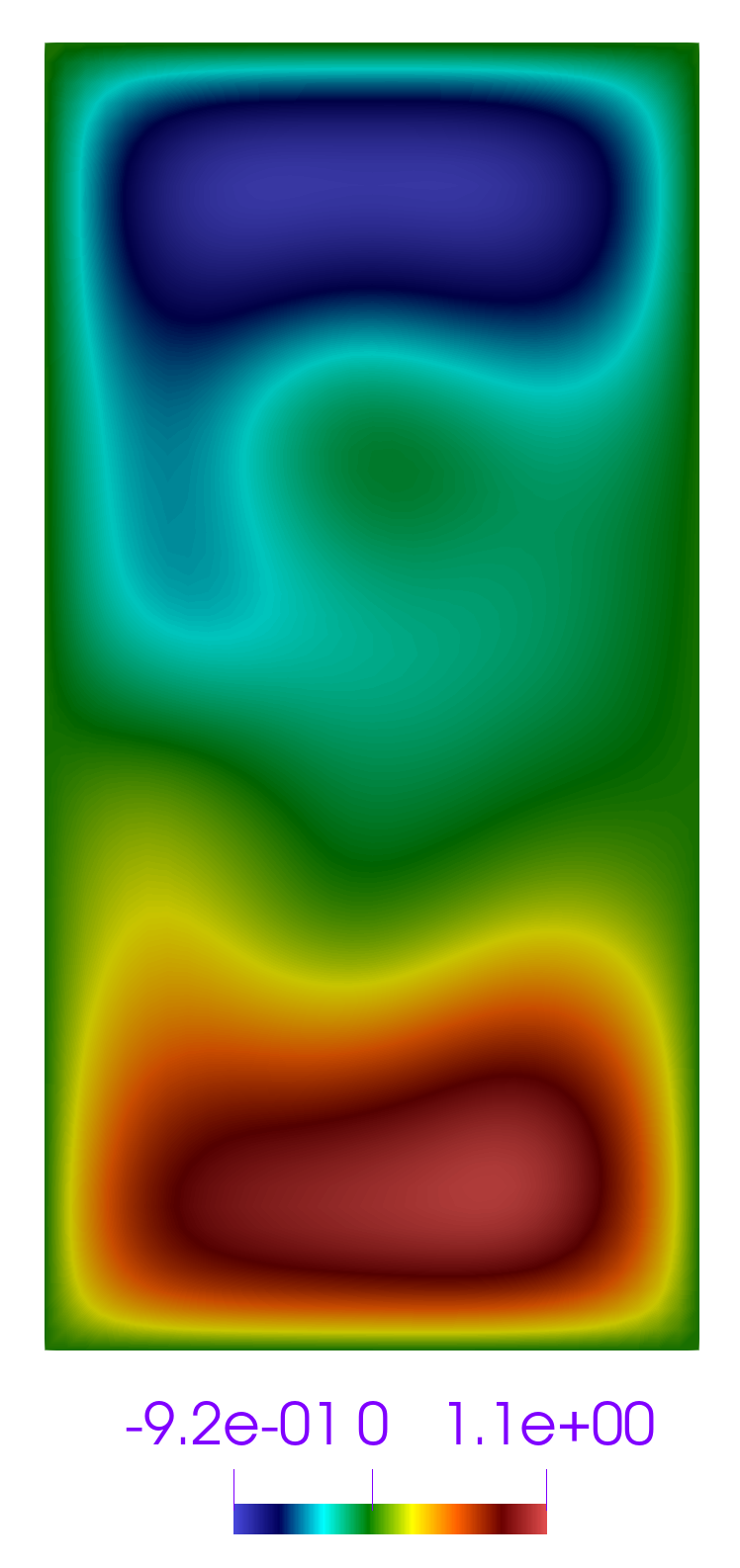}
    \end{overpic}
\put(-58,180){$N_q^r = 20$}
    \begin{overpic}[width=0.192\textwidth]{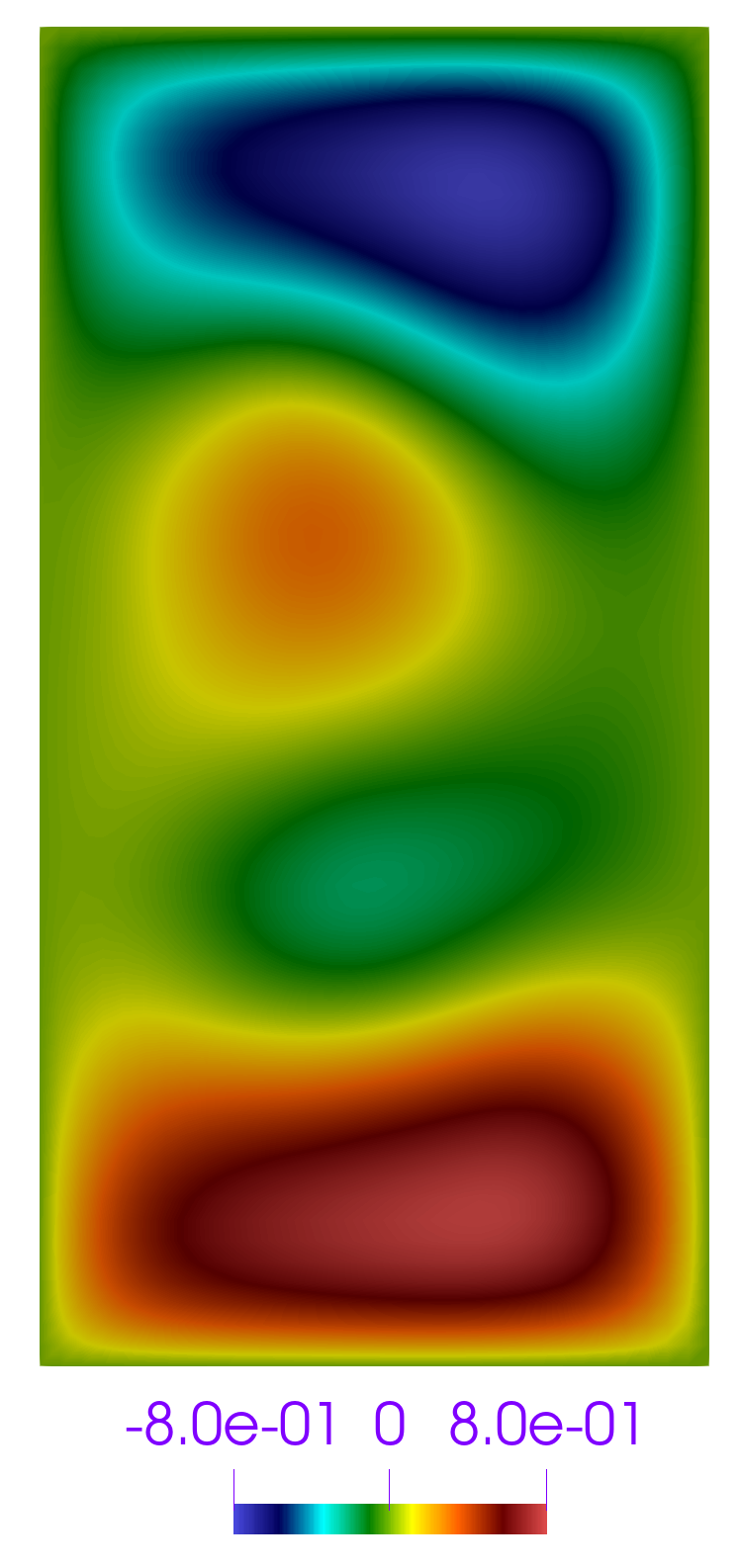}
    \end{overpic}
\put(-58,180){$N_q^r = 30$} \\
    \begin{overpic}[width=0.192\textwidth]{sections/img/FOM_case2.png}
    \end{overpic}
    \put(-58,180){FOM}
    \begin{overpic}[width=0.189\textwidth]{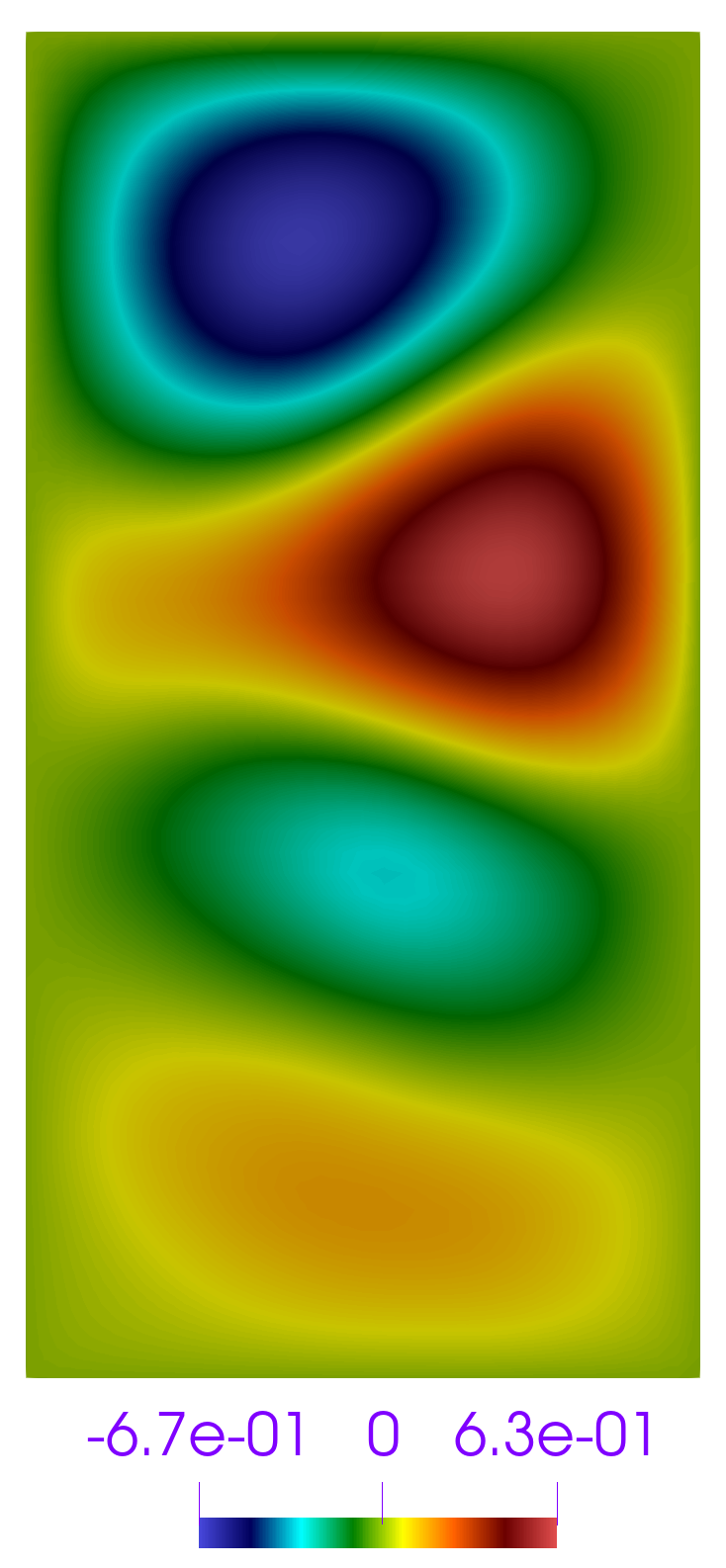}
    \end{overpic}
\put(-58,180){$N_q^r = 10$}
    \begin{overpic}[width=0.189\textwidth]{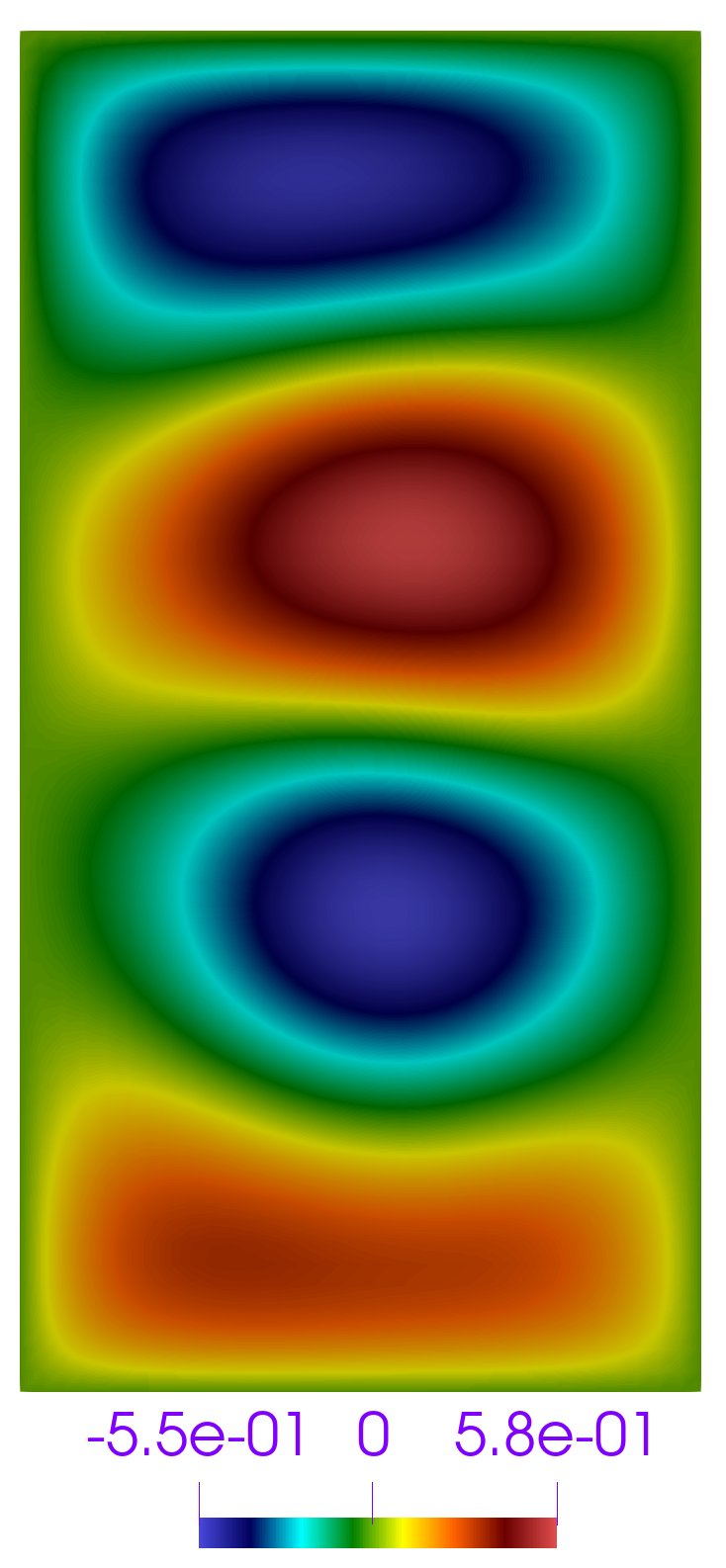}
    \end{overpic}
\put(-58,180){$N_q^r = 20$}
    \begin{overpic}[width=0.189\textwidth]{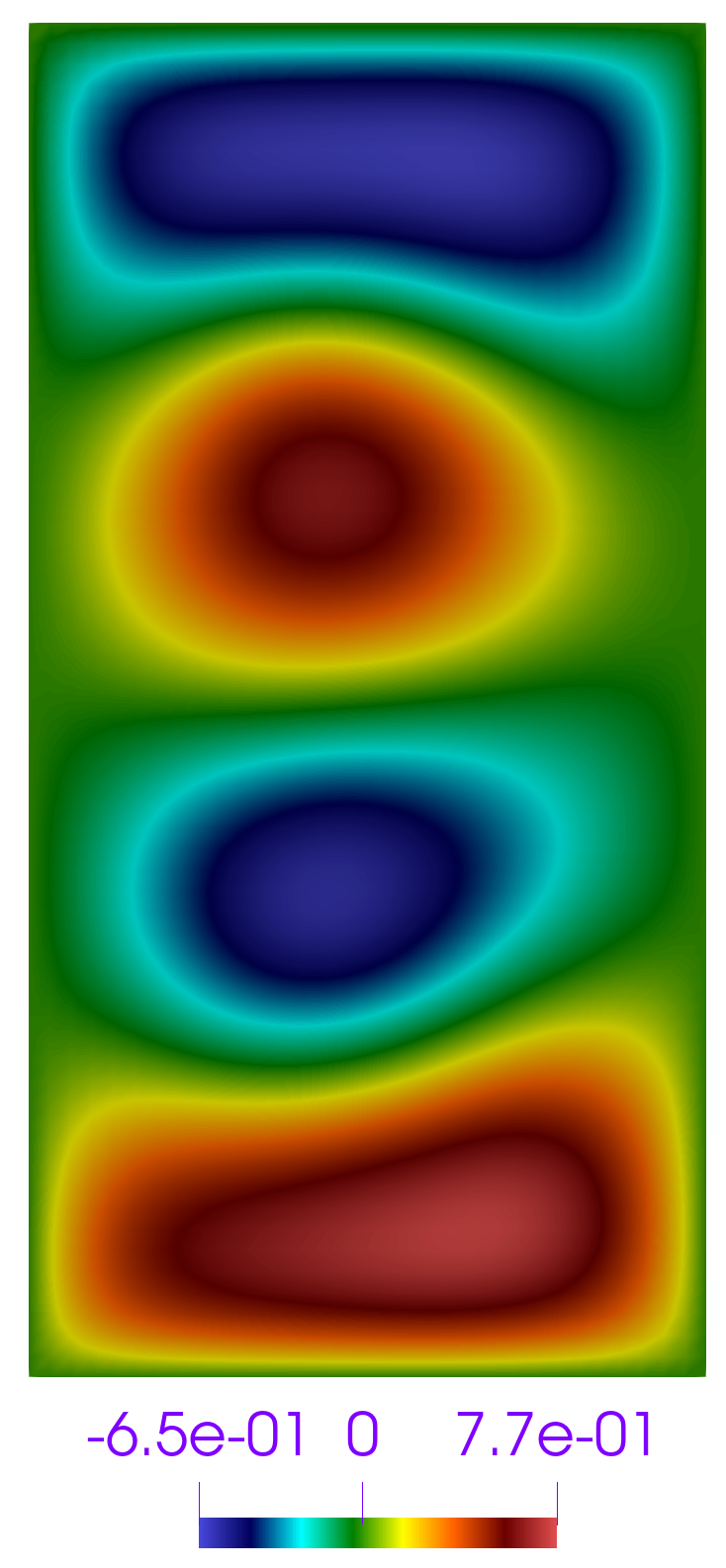}
    \end{overpic}
\put(-58,180){$N_q^r = 30$}
    \caption{
    Case 2: comparison of $\widetilde{\psi}$ computed by the FOM and the QGE-QGE ROM (top row)
    or the QGE-BV-$\alpha$ ROM (bottom row) for different numbers of POD vorticity modes $N_q^r$. $N_\psi^r$ is set to 10.} 
    \label{fig:psiROMstand2}
\end{figure}

Following what we have done for Case 1, we set 
$N^\psi_r = 10$ in order to retain 98\% of the eigenvalue energy associated with $\psi$ and we let $N^q_r$ vary. For Case 2, with $N^q_r$ = 10, 20, 30 we capture 59\%, 71\% and 76\% of the eigenvalue energy associated with $q$.
So, once again the ROM simulations are rather severely under-resolved.
Figure \ref{fig:energy_case2} compares the time evoution of the kinetic energy $E$ computed by the FOM and our two ROM approaches for the three values of $N_q^r$. There are substantial differences between Figure \ref{fig:energy_case2} and the corresponding figure for Case 1, i.e., Figure \ref{fig:energy_case1}. First of all, looking at the FOM kinetic energy, we see that for Case 2 it has oscillations with smaller amplitude and higher frequency (as one expects given the higher $Re$ in Case 2), while the average value is comparable in both cases. As for the ROMs, the
QGE-QGE ROM with $N_r^q = 10$ performs better in Case 2: while it provides a kinetic energy with larger amplitude and lower frequency than the FOM, the average is comparable to the ROM average. 
In Case 1 even the average of $E$ was off. 
On the other hand, the QGE-BV-$\alpha$ ROM with $N_q^r = 10$ seems to perform worse in Case 2, with several undershoots of the computed kinetic energy.
As $N_q^r$ is increased to 20 and 30, the kinetic energies computed by ROMs get closer to the FOM kinetic energy. This improvement is reflected in the computed $\widetilde{\psi}$ shown in Figure \ref{fig:psiROMstand2}.
Like in Case 1, we can see that the QGE-QGE ROM fails to capture the four-gyre pattern for $N_q^r = 10, 20$. However, when $N_q^r$ is increased to 30, such pattern starts to emerge. We note that that $\widetilde{\psi}$ computed by the QGE-QGE ROM with $N_q^r = 30$ for Case 2 (last panel in the first row of Figure \ref{fig:psiROMstand2}) looks similar to the $\widetilde{\psi}$ computed by the same ROM with $N_q^r = 40$ for Case 1 (second panel in Figure \ref{fig:psiROMalpha})
This is not surprising since the retained eigenvalue energy is comparable: 76\% for the former vs 74\% for the latter. 
From the panels on the bottom row of Figure 
\ref{fig:psiROMstand2}, we see that the pattern in $\widetilde{\psi}$ given by the QGE-BV-$\alpha$ ROM matches the FOM pattern already for $N_q^r = 20$. When $N_q^r$ is increased, the QGE-BV-$\alpha$ ROM becomes less diffusive and the magnitude of $\widetilde{\psi}$ computed by the ROM gets closer to the magnitude computed by the FOM. 

\begin{figure}[htb]
        \begin{overpic}[width=0.19\textwidth]{sections/img/FOM_case2.png}
    \end{overpic}
    \put(-58,180){FOM}
        \begin{overpic}[width=0.19\textwidth]{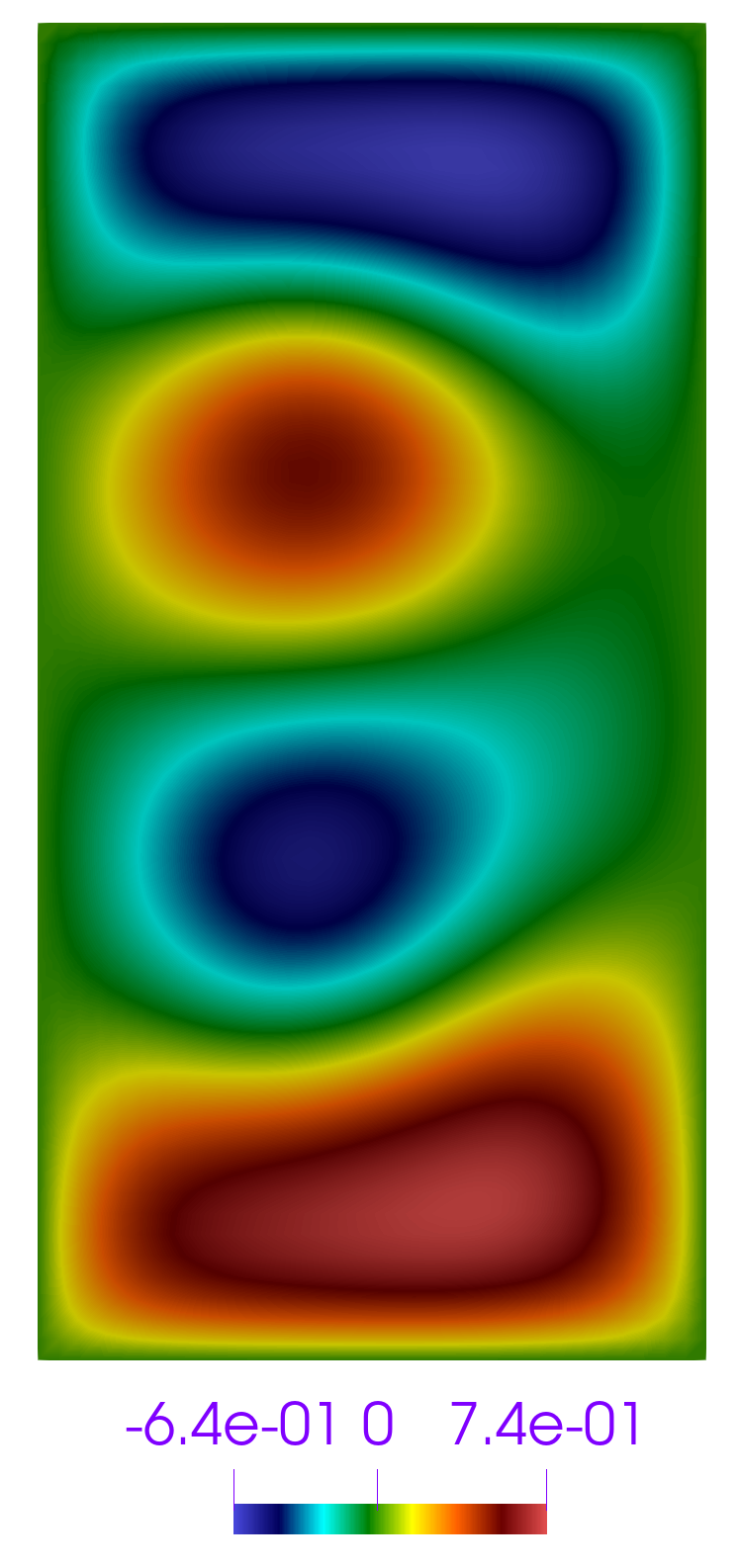}
    \end{overpic}
\put(-58,180){$N_q^r = 40$}
    \caption{Case 2: $\widetilde{\psi}$ computed by the FOM and the QGE-QGE ROM with $N_q^r = 40$ and $N_\psi^r = 10$.}
    \label{fig:psiROMalpha2}
\end{figure}

For a more quantitative comparison, we report in Table \ref{tab:errors2} the $L^2$ error \eqref{eq:normerror} for both ROMs with different values of $N_q^r$. Like in Case 1, we see that the errors obtained with the QGE-BV-$\alpha$ ROM are smaller in all cases, although the errors of the QGE-QGE ROM and QGE-BV-$\alpha$ ROM get closer as $N_q^r$ increases. We note that if $N_q^r$ is set to 40 for the QGE-QGE ROM (corresponding to retaining  80\% of the eigenvalue energy), the computed $\widetilde{\psi}$ shown in Figure \ref{fig:psiROMalpha2} gets very close
to the QGE-BV-$\alpha$ solution with $N_q^r = 30$.

\begin{table}[htb]
\centering
\begin{tabular}{|c|c|c|c|}
\multicolumn{4}{c}{} \\
\cline{1-4}
$N_q^r$ & Relative energy content & $\varepsilon_{QGE-QGE}$ & $\varepsilon_{QGE-BV-\alpha}$\\
\hline
10 & 59\% & 5e+00 & 7.4e-01 \\
20 & 71\% & 1e+00 & 5.3e-01 \\
30 & 76\% & 6.4e-01 & 3.7e-01\\
\hline
\end{tabular}
\caption{
Case 2: $L^2$ error \eqref{eq:normerror} given by QGE-QGE-ROM and QGE-BV-$\alpha$ ROM for $N_q^r = 10, 20, 30$ and $N_\psi = 10$.}
\label{tab:errors2}
\end{table}

We conclude this section with a comment on the efficiency of our ROM approaches in the spirit of Table \ref{tab:time}. Table \ref{tab:time2} reports the CPU time required by the ROMs in Case 2 and the relative speed up with respect to the FOM. 
First of all, we notice that the CPU times required by the ROMs are very similar for Case 1 and 2, while 
we observe more important speed-ups for Case 2. This is consistent with our hypothesis at the end of Sec.~\ref{sec:case1}, i.e., that the use of a finer mesh would lead to larger computational savings. In fact, we recall that for Case 2 we use mesh 32$\times$64, while for Case 1 we took mesh 16$\times$32. 
Furthermore, it explains why previous work that used a 256 $\times$ 512 mesh \cite{San2014} reports a speed-up of the order of 100.  

\begin{table}[htb]
\centering
\begin{tabular}{|c|c|c|}
\hline
& $\varepsilon_{QGE-QGE}$ & $\varepsilon_{QGE-BV-\alpha}$\\
\hline
CPU time & 113 s & 177 s \\
\hline
speed-up & 10.3 & 6.6\\
\hline
\end{tabular}
\caption{Case 2: CPU time required by 
QGE-QGE-ROM and QGE-BV-$\alpha$ ROM with $N_q^r = 10$ and $N_\psi = 10$ and relative speed-up with respect to the CPU time required by the FOM simulation (1166 s).}
\label{tab:time2}
\end{table}

\section{Concluding remarks}\label{sec:End}
This paper presents a novel regularization for reduced order models of the quasi-geostrophic equations when the POD modes retained to construct the reduced basis are insufficient to describe the system dynamics. The proposed regularization draws inspiration from the linear BV-$\alpha$ model, which has been used only as a replacement of the QGE at the full order level so far. For the collection of the snapshots, we apply a Finite Volume method, which has the advantage of preserving conservation of conserved quantities at the discrete level.
The particular ROM approach that we combined with the new regularization is of the POD-Galerkin type. To show the effectiveness of the BV-$\alpha$ closure model, we compare the results computed by the ROM with and without regularization for the classical double-gyre wind forcing benchmark. We consider two cases with the same Munk scale, one with small Rossby number and the other with high Reynolds number. 

Our numerical results show that for both cases the solution computed by the regularized ROM is more accurate, even when the retained POD modes account for a small percentage of the cumulative eigenvalue energy (i.e., about 50-60\%). The price to pay for this increased accuracy is an increased computational cost: the CPU time of the regularized ROM is about 1.5 times the CPU time required by the non-regularized ROM. Despite this increased computational cost, the regularized ROM is still a competitive alternative to the full order model. In fact, its cost is about 1/3 (for the small $Ro$ case) to 1/6 (for the high $Re$ case) of the cost of the full order model.

While the ROM regularized by the linear BV-$\alpha$ model is accurate in the reconstruction of the solution patterns, the positive and negative peaks in the magnitude get smoothed out. This is expected since linear filters are known to be dissipative. At the full order level, we have shown the a nonlinear version of the BV-$\alpha$ model introduces much less artificial dissipation \cite{Girfoglio1}. Hence, a natural extension of the work presented in this paper is a regularization inspired by this nonlinear BV-$\alpha$ model.


\section*{Acknowledgements}
We acknowledge the support provided by the European Research Council Executive Agency by the Consolidator Grant project AROMA-CFD ``Advanced Reduced Order Methods with Applications in Computational Fluid Dynamics" - GA 681447, H2020-ERC CoG 2015 AROMA-CFD, PI G. Rozza, ARIA H2020-MSCA-RISE-2019 project, PON ``Research and Innovation on Green related issues" FSE REACT-EU 2021 project, PRIN NA\_FROM-PDEs project and INdAM-GNCS 2019-2021 projects.
This work was also partially supported by US National Science Foundation through grant DMS-1953535 (PI A.~Quaini). 
A.~Quaini also acknowledges support from the Radcliffe Institute for Advanced Study at Harvard University where
she has been a 2021-2022 William and Flora Hewlett Foundation Fellow.





\bibliographystyle{plain}
\bibliography{sample}

\end{document}